\documentclass[final]{siamltex}
\usepackage{graphicx} % Required for inserting images
\usepackage{mathtools}
\usepackage{latexsym}
\usepackage[mathscr]{euscript}
\usepackage{xcolor}
\usepackage{algorithm}
\usepackage{fancyhdr}
\usepackage{multirow}
\usepackage{epsfig}
\usepackage{mathptmx}       % selects Times Roman as 
\usepackage[most]{tcolorbox}
\usepackage{helvet}         % selects Helvetica as sans-serif font
\usepackage{courier}        % selects Courier as typewriter font
\usepackage{type1cm}        % activate if the above 3 fonts are
\usepackage{imakeidx}         % allows index generation
\usepackage[a4paper,textheight=23cm,textwidth=15cm]{geometry}
\usepackage{algpseudocode}
\usepackage{geometry}
\usepackage{hyperref}
\usepackage{enumitem}
\usepackage{titlesec}
\usepackage{amsmath,amssymb,amsfonts}
\usepackage{hyperref}
\usepackage{tikz}
\usepackage{bm}
\usepackage{booktabs}
\newtheorem{remark}{remark}
\newcommand{\R}{\mathbb{R}}

\newcommand{\nd}{\noindent}

\title{A survey  of scalar and vector extrapolation }
\author{K. Jbilou\thanks{Université du Littoral Cote d'Opale, LMPA, 50 rue F. Buisson, 62228 Calais-Cedex, France.} }
\date{}

\begin{document}

\maketitle

\begin{abstract}
Scalar extrapolation and convergence acceleration methods are central tools in numerical analysis for improving the efficiency of iterative algorithms and the summation of slowly convergent series. These methods construct transformed sequences that converge more rapidly to the same limit without altering the underlying iterative process, thereby reducing computational cost and enhancing numerical accuracy. Historically, the origins of such techniques can be traced back to classical algebraic methods by Al-Khwarizmi and early series acceleration techniques by Newton, while systematic approaches emerged in the 20th century with Aitken's $\Delta^2$ process and Richardson extrapolation. Later developments, including the Shanks transformation and Wynn’s $\varepsilon$-algorithm, provided general frameworks capable of eliminating multiple dominant error components, with deep connections to Padé approximants and rational approximations of generating functions. This paper presents a comprehensive review of classical scalar extrapolation methods, including Richardson extrapolation, Aitken’s $\Delta^2$ process, Shanks transformation, Wynn’s $\varepsilon$-algorithm, and $\rho$ and $\theta$ algorithms.  We examine their theoretical foundations, asymptotic error models, convergence properties, numerical stability, and practical implementation considerations. The second part of this work is dedicated to vector extrapolation methods: polynomial based ones and $\varepsilon$-algorithm generalizations to vector sequences. 
Additionally, we highlight modern developments such as  their applications to iterative solvers, Krylov subspace methods, and large-scale computational simulations. The aim of this review is to provide a unified perspective on scalar and vector extrapolation techniques, bridging historical origins, theoretical insights, and contemporary computational applications.
\end{abstract}

{\bf keywords}: Scalar extrapolation, Convergence acceleration, Sequence transformations, Aitken's $\Delta^2$ process, Richardson extrapolation, Shanks transformation,  $\varepsilon$-algorithm, Padé approximants, Vector and matrix extrapolation.

\section{Historical background of  extrapolation methods}
Scalar extrapolation methods, also referred to as sequence transformations or convergence acceleration techniques, are fundamental tools in numerical analysis. Their primary objective is to accelerate the convergence of sequences $\{x_n\}$ toward a limit $x^\ast$, or to assign finite values to slowly convergent or divergent series. These methods have broad applications in numerical solutions of differential equations, iterative linear and nonlinear solvers, series summation, and approximation theory \cite{Brezinski1991, Sidi2003}.

%\subsection{Early Origins and Motivations}

The concept of iterative approximation is ancient. In the 9th century, Al-Khwarizmi \cite{al-Khwarizmi} introduced systematic methods for solving algebraic equations. His algorithmic approach can be viewed as a precursor to modern iterative methods, emphasizing step-by-step computation toward an approximate solution. Later, mathematicians such as Newton (1642–1727) and James Gregory (1638–1675) developed series expansions, interpolation, and iterative formulas for finding roots of equations, which implicitly relied on accelerating convergence to a solution.

These early methods were largely motivated by practical needs: engineers, astronomers, and physicists required accurate approximations of roots, integrals, and series sums without access to modern computers. The challenge was to reduce computational effort while achieving sufficient precision, laying the foundation for the formal study of sequence transformations.
\\
The formalization of convergence acceleration began in the 20th century: Aitken's $\Delta^2$ Process (1926).  Aitken \cite{Aitken1926} introduced a simple formula to accelerate linearly convergent sequences. By eliminating the leading term of the error, his method significantly reduces the number of iterations required to achieve a desired accuracy. This process is particularly effective for sequences whose errors behave approximately like a single exponential term.

In 1911,  Richardson \cite{Richardson1911} proposed a systematic approach for sequences depending on a parameter $h$, often representing step size in numerical differentiation or integration. By assuming an asymptotic expansion of the form $x(h) = x^\ast + c_1 h + c_2 h^2 + \dots$, Richardson's method constructs linear combinations of sequence terms at different $h$ values to cancel leading-order errors. This idea is foundational in high-accuracy numerical methods and remains widely used in modern computational mathematics. 
Both Aitken and Richardson addressed the practical problem of reducing computational effort by intelligently predicting the limit of a sequence without requiring large numbers of iterations.

New transformations and Wynn's $\varepsilon$-algorithm
appear during the Mid-20th Century. 
In the 1950s, more sophisticated techniques were developed to handle sequences with multiple error components. In 1955,  Shanks \cite{Shanks1955} generalized Aitken’s method to sequences whose errors contain several exponential terms. The transformation can be expressed in a determinant form, making it theoretically elegant but computationally expensive. Shanks showed that his transformation could significantly improve convergence in a variety of numerical problems.
    
    In 1956),  Wynn \cite{Wynn1956} introduced a recursive algorithm called $\varepsilon$-algorithm, that implements the Shanks transformation efficiently. The $\varepsilon$-algorithm avoids direct computation of determinants, making it numerically stable and suitable for practical applications. This algorithm can eliminate multiple error terms iteratively, providing accelerated convergence for sequences arising in iterative solvers, series summation, and other numerical computations. For more details see \cite{Brezinski1971,Brezinski1971r,Brezinski1980,Brezinski1991,Wimp1981,Wynn1962,
    Wynn1966}.

One of the most profound insights in the theory of sequence transformations is their connection to Padé approximants. The determinant forms of Shanks transformations and the recursive formulas of Wynn can be interpreted as constructing rational approximations to the generating function of a sequence \cite{Baker1981,Brezinski1980,Brezinski1983}. This connection explains why these transformations can effectively sum series that converge slowly or even diverge, especially when the series represents the expansion of a function near a singularity.
\\
The historical evolution of scalar extrapolation methods reflects both the practical needs and the mathematical sophistication of numerical analysis:

\begin{itemize}
    \item Early methods focused on simple iterative improvement and practical computation.  
    \item 20th-century methods formalized the elimination of leading error terms and established systematic procedures for convergence acceleration.  
    \item Modern developments extend these techniques to high-dimensional, vector, and matrix sequences, maintaining both efficiency and numerical stability.
\end{itemize}

Understanding this evolution provides crucial insight into how classical methods inform modern algorithms and why sequence transformations remain a fundamental tool in numerical analysis. 
Many numerical procedures generate sequences that converge slowly to a desired limit. Such sequences commonly arise in fixed-point iterations, numerical quadrature, discretizations of differential equations, and series summation. Improving convergence rates is essential in reducing computational cost and increasing numerical efficiency.
\\
Scalar extrapolation methods provide a non-intrusive approach to convergence acceleration by post-processing the generated sequence. Their simplicity and effectiveness have made them classical yet enduring tools in numerical analysis.

While the classical methods discussed above primarily target scalar sequences, many modern computational problems involve vector or matrix sequences, such as those arising from iterative solvers of linear and nonlinear systems, or from discretized PDEs. Scalar extrapolation techniques, including Aitken’s $\Delta^2$ process, Richardson extrapolation, and Wynn’s $\varepsilon$-algorithm, provide the theoretical foundation for their vector counterparts. Methods such as the vector and topological $\varepsilon$-algorithm, minimal polynomial extrapolation (MPE), and reduced rank extrapolation (RRE) extend the principles of sequence transformations to higher-dimensional settings \cite{Brezinski1991,Jbilousadok2000,Sidi2003}. These vector and matrix extrapolation techniques are designed to accelerate convergence of sequences of vectors or matrices by eliminating dominant error components in a similar manner to scalar methods, while accounting for linear dependence and multi-dimensional correlations. They are particularly valuable in large-scale simulations and iterative solvers, where each iteration is computationally expensive, and rapid convergence can significantly reduce overall computational cost. Understanding scalar extrapolation is thus essential, as it provides both the conceptual framework and analytical tools for these higher-dimensional generalizations.

The remainder of this paper is organized as follows.  
Section~2 reviews several well-known scalar extrapolation methods, including
Richardson extrapolation, Aitken’s $\Delta^2$ process, the Shanks transformation,
Wynn’s scalar $\varepsilon$-algorithm, as well as the $E$-, $\rho$-, and
$\theta$-algorithms.  
Section~3 is devoted to classical vector extrapolation techniques. In particular,
we present the Reduced Rank Extrapolation (RRE), the Minimal Polynomial
Extrapolation (MPE), the Modified Minimal Polynomial Extrapolation (MMPE), the
vector and topological $\varepsilon$-algorithms (VEA and TEA), and Anderson
acceleration (AA). We also present some applications of these vector extrapolation methods for linear and nonlinear problems. 
Finally, Section~4 illustrates the performance of these methods through several
numerical examples.

\section{Scalar extrapolation methods}

\subsection{Richardson extrapolation}

In numerical analysis, \emph{Richardson extrapolation} \cite{Richardson1911, Richardson1927} is a classical technique for accelerating the convergence of sequences or iterative methods. 
Convergence acceleration refers to a collection of methods aimed at improving the rate at which a sequence or an iterative approximation approaches its exact limit or solution. When an iterative method produces a sequence of successive approximations, acceleration techniques transform this sequence to produce a more accurate estimate with fewer iterations, thereby reducing computational cost. Richardson extrapolation is a prototypical example: it combines several approximations computed using different step sizes to eliminate leading-order error terms, thereby enhancing the precision of the final estimate.
\\
This method plays a central role in many areas of computational science. It is widely used in numerical integration, the numerical solution of differential equations, and physical simulations, where computational efficiency is critical. By reducing the number of iterations or the step size required for a given accuracy, Richardson extrapolation can also help mitigate numerical errors arising from finite-precision arithmetic or, in future contexts, quantum computing.
A particularly important application of Richardson extrapolation is in the construction of high-accuracy numerical integration schemes. For instance, Romberg integration is essentially an accelerated version of the trapezoidal rule obtained by successive Richardson extrapolations, which systematically eliminate leading-order discretization errors and achieve rapid convergence to the exact integral.

\medskip 
\nd Assume an asymptotic expansion of the form
\begin{equation}
x(h) = x^\ast + c_1 h^p + c_2 h^{p+1} + c_3 h^{p+2} + \cdots,
\label{eq:richardson_asymptotic}
\end{equation}
where $p>0$ and $h \to 0$.
The key idea is to eliminate the leading error term $c_1 h^p$ using approximations at two different step sizes.

%\subsection{Two-Point Richardson Extrapolation}

\nd Let $x(h)$ and $x(h/2)$ be approximations computed at step sizes $h$ and $h/2$. Using \eqref{eq:richardson_asymptotic}, we have
\begin{align*}
x(h) &= x^\ast + c_1 h^p + c_2 h^{p+1} + \cdots, \\
x(h/2) &= x^\ast + c_1 \left(\frac{h}{2}\right)^p + c_2 \left(\frac{h}{2}\right)^{p+1} + \cdots.
\end{align*}
Solving for $x^\ast$ by eliminating $c_1$ gives the Richardson extrapolated value:
\begin{equation}
x^* \approx \frac{2^p x(h/2) - x(h)}{2^p - 1}.
\label{eq:richardson}
\end{equation}
This formula cancels the $h^p$ term, resulting in an approximation whose leading error term is now $\mathcal{O}(h^{p+1})$.
\\
Richardson extrapolation can be applied recursively to obtain higher-order convergence. Let
\[
t_0^{(n)} = x(h/2^n),
\]
and define recursively for $k \ge 1$:
\begin{equation}
t_k^{(n)} = t_{k-1}^{(n+1)} + \frac{t_{k-1}^{(n+1)} - t_{k-1}^{(n)}}{2^{p+k-1} - 1}.
\end{equation}
Then $t_k^{(0)}=s_n$ approximates $x^\ast$ with error $\mathcal{O}(h^{p+k})$.
\\
More generally, Richardson extrapolation produces a sequence $(t_k^{(n)})$ difined as follows
\begin{equation}\label{rich2}
t_k^{(n)}=t_{k-1}^{(n)} - \displaystyle \frac{t_{k-1}^{(n+1)}-t_{k-1}^{(n)}}{x_{n+k}-x_n}\,x_n,\;\; t_0^{(n)}=s_n\;  n\ge 0,\; k\ge 1,
\end{equation}
where $(x_n)$, $n\ge 0$ are  given points. The Richardson process consists in interpolating the sequence terms $s_n,\ldots,s_{n+k}$ at the points $x_n,\ldots,x_{n+k}$. 
% Assume the asymptotic expansion
% \[
% x(h) = x^\ast + \sum_{j=1}^{\infty} c_j h^{p+j-1}.
% \]
% The first-order Richardson extrapolation reduces the error from
% \[
% e_1 = c_1 h^p + \mathcal{O}(h^{p+1})
% \]
% to
% \[
% e_2 = \frac{c_2}{2^p - 1} h^{p+1} + \mathcal{O}(h^{p+2}).
% \]
% By recursively applying the procedure, one can systematically eliminate higher-order terms and achieve arbitrarily high-order accuracy, provided the expansion exists and $h$ is sufficiently small.
\\
A classical application is Romberg integration. Let $I(h)$ denote a trapezoidal approximation of an integral with step size $h$:
\[
I = \int_a^b f(x)\,dx.
\]
Using Richardson extrapolation with $I(h)$ and $I(h/2)$, one obtains
\[
I \approx  \frac{4 I(h/2) - I(h)}{3},
\]
which improves the convergence from $\mathcal{O}(h^2)$ to $\mathcal{O}(h^4)$. 
Richardson extrapolation is a robust and widely used technique to improve the convergence of sequences with algebraic error terms. Its applications range from numerical integration and differentiation to the solution of differential equations and iterative methods. The method is simple, systematic, and can achieve arbitrarily high-order accuracy under appropriate conditions.\\
Convergence acceleration results of the Richardson sequence $(t_k^{(n)})$ are well known \cite{Laurent1964}. The main theorem says that a necessary and sufficient condition that $\displaystyle \lim_{n \ \leftarrow \infty} t_k^{(n)}=s$ for all sequences $(s_n)$ that converge to $s$ is that there exist $\alpha$ and $\beta$ with $\alpha <1< \beta$ such that $\forall n,\; x_{n+1}/x_n \notin [\alpha,\, \beta]$.

\subsection{Aitken’s $\Delta^2$ process}

\subsubsection{The classical $\Delta^2$ algorithm}

Aitken’s $\Delta^2$ process is a classical scalar extrapolation technique designed to accelerate the convergence of a linearly convergent sequence. Let $\{s_n\}$ be a sequence converging to a limit $x^\ast$ with linear convergence rate $\rho$, i.e.,
\begin{equation}
\lim_{n \to \infty} \frac{|s_{n+1}-s|}{|s_n-s|} = \rho, \quad 0 < \rho < 1.
\end{equation}
Define the forward differences:
\begin{align}
\Delta s_n &= s_{n+1} - s_n, \\
\Delta^2 s_n &= s_{n+2} - 2 s_{n+1} + s_n.
\end{align}
The Aitken transformed sequence is given by 
\begin{equation}
t_n = s_n - \frac{(\Delta s_n)^2}{\Delta^2 s_n}, \quad \Delta^2 s_n \neq 0.
\label{eq:aitken}
\end{equation}
Intuitively, the method extrapolates the linear trend in consecutive sequence elements to estimate the limit more accurately.
\\
% Suppose $\{x_n\}$ admits a linear error model
% \begin{equation}
% x_n - x^\ast = \rho^n (x_0 - x^\ast) + \epsilon_n, \quad |\epsilon_n| \ll |x_0 - x^\ast|.
% \end{equation}
% Then
% \begin{align*}
% \Delta x_n &= x_{n+1} - x_n = (\rho^{n+1} - \rho^n)(x_0 - x^\ast) + (\epsilon_{n+1} - \epsilon_n) \\
% &= \rho^n (\rho - 1)(x_0 - x^\ast) + \mathcal{O}(\epsilon_n), \\
% \Delta^2 x_n &= x_{n+2} - 2x_{n+1} + x_n = \rho^n (\rho - 1)^2 (x_0 - x^\ast) + \mathcal{O}(\epsilon_n).
% \end{align*}
% Substituting into \eqref{eq:aitken} gives the accelerated error:
% \begin{equation}
% \hat{x}_n - x^\ast = \mathcal{O}(\epsilon_n).
% \end{equation}
\emph{Kernel of a sequence transformation.}
Let $\mathcal{T}$ be a sequence transformation that associates with a scalar
sequence $(s_n)$ a transformed sequence $t_n=\mathcal{T}(s_n)$, and assume that
$(s_n)$ converges to a limit $s$.

\medskip
\noindent
\textbf{Definition (Kernel).}
The \emph{kernel} of the transformation $\mathcal{T}$, denoted by
$\mathcal{K}(\mathcal{T})$, is defined as
\[
\mathcal{K}(\mathcal{T})
=
\left\{
(s_n)\;:\; t_n = s
\text{ for all } n \text{ for which the transformation is defined}
\right\}.
\]
Equivalently, introducing the error sequence $e_n = s_n - s$, the kernel can be
characterized as
\[
\mathcal{K}(\mathcal{T})
=
\left\{
(e_n)\;:\; (\mathcal{T}(e_n) = 0
\right\},
\]
that is, the set of all error sequences annihilated by the transformation. For the Aitken--$\Delta^2$ transformation,  the kernel is given by
\[
\mathcal{K}=
\left\{
(s_n)\;:\;
s_n = s + c\,\lambda^{n}, \qquad
c \neq 0,\;\; \lambda \neq 1
\right\}.
\]

\medskip
\nd The application of Aitken’s $\Delta^2$ process presupposes that the sequence
$\{s_n\}$ converges to a limit $s$ and that its error admits an asymptotic
representation of the form
\[
e_n = s_n - s = c\,\lambda^{\,n} + o(\lambda^{\,n}),
\]
where $c \neq 0$ and the convergence factor satisfies $0 < |\lambda| < 1$.
Under this hypothesis, the original sequence is linearly convergent. 
Within this setting, the $\Delta^2$ transformation produces a new sequence
$(t_n)$ that converges to the same limit $s$, but with an improved rate of
convergence. In particular, if the ratio of successive errors
\[
\frac{e_{n+1}}{e_n}
\]
tends to a limit $\lambda$ with $|\lambda|<1$, the leading term of the error
expansion is effectively cancelled by the transformation. Thus, in the absence of higher-order terms, the Aitken process effectively eliminates the leading linear error term, transforming linear convergence into quadratic convergence; see \cite{Brezinski1971} for more details on Atken's process.

% Notice that Richardson extrapolation can be seen as a special case of Aitken's $\Delta^2$ process  applied to sequences whose errors behave as powers of $h$. Moreover, applying Richardson extrapolation recursively corresponds to constructing Padé-type rational approximations in the discretization parameter $h$, providing a theoretical justification for its convergence enhancement.
% \subsection{Generalization to Nonlinear Convergence}
% \\
% For sequences with errors of the form
% \begin{equation}
% x_n - x^\ast = c_1 \lambda_1^n + c_2 \lambda_2^n + \cdots, \quad |\lambda_1| > |\lambda_2| > \cdots,
% \end{equation}
% the Aitken $\Delta^2$ method removes the dominant term $c_1 \lambda_1^n$, improving convergence toward $x^\ast$. However, the presence of multiple significant exponential terms may require repeated or higher-order transformations, such as the Shanks transformation. \\
% Applying Aitken’s $\Delta^2$ process recursively leads to higher-order acceleration. Let
% \[
% \hat{x}_n^{(0)} = x_n, \quad 
% \hat{x}_n^{(k+1)} = \hat{x}_n^{(k)} - \frac{\left(\Delta \hat{x}_n^{(k)}\right)^2}{\Delta^2 \hat{x}_n^{(k)}}.
% \]
% If the sequence satisfies a linear or exponential model, iterated Aitken transformations can eliminate successive dominant terms, further enhancing convergence.
\noindent Geometrically, Aitken’s method can be interpreted as projecting the current sequence along a line defined by consecutive points to estimate the intersection with the limit. This viewpoint explains the method’s effectiveness for sequences exhibiting nearlinear behavior in consecutive iterates.
\\
%\subsection{Numerical Stability Considerations}
Although Aitken’s $\Delta^2$ process is simple, numerical instability may occur when
\[
\Delta^2 x_n = x_{n+2} - 2x_{n+1} + x_n \approx 0,
\]
leading to large round-off errors. Remedies include:
\begin{itemize}
    \item Avoiding acceleration if $\Delta^2 x_n$ is below a tolerance,
    \item Combining with smoothing or averaging techniques,
    \item Using the Shanks transformation or $\varepsilon$-algorithm for sequences with multiple exponential components.
\end{itemize}
%\subsection{Summary of Convergence Improvement}
% For linearly convergent sequences:
% \[
% |x_n - x^\ast| = \mathcal{O}(\rho^n) \quad \Rightarrow \quad |\hat{x}_n - x^\ast| = \mathcal{O}(\rho^{2n}).
% \]
% Hence, Aitken’s $\Delta^2$ process often provides dramatic reduction in iteration count in practical numerical problems.

\subsubsection{The iterated \(\Delta^2\) method}
The \emph{iterated \(\Delta^2\) method} is a classical sequence transformation
introduced by Aitken for accelerating the convergence of sequences
\(\{x_n\}_{n\ge 0}\), particularly those that converge linearly. Its main idea
is to eliminate the leading term of the asymptotic error, producing a sequence
that converges faster to the limit.
\medskip
\noindent
The efficiency of the basic \(\Delta^2\) method can be improved by applying it
iteratively. Denote the \(k\)-th iterated sequence by \(t_k^{(n)}\), with
\begin{equation}
t_0^{(n)} = s_n, \qquad 
t_{k+1}^{(n)} = t_k^{(n)} - \frac{(\Delta t_k^{(n)})^2}{\Delta^2 t_k^{(n)}}, \qquad k \ge 0.
\label{eq:iterated_delta2}
\end{equation}
Each iteration attempts to remove the leading term of the remaining error in
the transformed sequence. In practice, a few iterations often suffice to achieve
a substantial acceleration, especially for sequences with approximately linear
convergence.

\medskip
\noindent 
 The iterated \(\Delta^2\) method is particularly effective for
    sequences with errors that can be modeled as \(e_n \approx C \lambda^n\)
    with \(|\lambda|<1\), i.e., linearly convergent sequences. It is a scalar method but can be extended to vector sequences,
    yielding a simple form of vector extrapolation.
The method is widely used for accelerating the convergence of
    iterative methods for nonlinear equations, summation of series, and
    sequence extrapolation in numerical analysis.

\medskip
\noindent  
The iterated \(\Delta^2\) method is the simplest form of sequence
transformation and can be seen as a precursor to more general nonlinear
accelerators such as Wynn’s \(\varepsilon\)-algorithm, the \(\theta\)-algorithm,
and Anderson acceleration. In particular, the \(\varepsilon\)-algorithm can be
interpreted as a systematic generalization of the iterated \(\Delta^2\) process
to arbitrary orders of extrapolation; this will be seen in the coming sections.

\subsection{Shanks transformation}
The classical Shanks transformation \cite{Shanks1955} accelerates convergence of a scalar sequence $\{x_n\}$ by eliminating the leading exponential term in the error. The \emph{general Shanks transformation} extends this idea to remove multiple dominant error terms simultaneously, allowing acceleration of sequences with multiple exponential components:
\begin{equation}
s_n = s^\ast + \sum_{j=1}^m c_j \lambda_j^n, \quad |\lambda_1| > |\lambda_2| > \cdots > |\lambda_m|.
\end{equation}

\nd Let $\{s_n\}$ be a scalar sequence. The $k$-th order Shanks transformation, denoted by $e_n^{(k)}$, aims to remove the first $k$ exponential components of the error. For $k=1$, it reduces to the classical Aitken $\Delta^2$ process:
\begin{equation}
e_n^{(1)} = s_n - \frac{(\Delta s_n)^2}{\Delta^2 x_n}, \quad \Delta s_n = s_{n+1}-x_n.
\end{equation}

\nd Higher-order transformations can be defined recursively using Wynn's $\varepsilon$-algorithm or via determinant expressions.
\\
The general Shanks transformation of order $k$ can be expressed as a ratio of two determinants. Let $\Delta^j x_n$ denote the $j$-th forward difference:
\[
\Delta^0 s_n = s_n, \quad \Delta^j s_n = \Delta^{j-1} s_{n+1} - \Delta^{j-1} s_n, \quad j \ge 1.
\]
Then the $k$-th order Shanks transform $e_k^{(n)}$ is given by
\begin{equation}
e_k^{(n)} = e_k(s_n)= \frac{
\begin{vmatrix}
s_n & s_{n+1} & \cdots & s_{n+k} \\
\Delta s_n & \Delta s_{n+1} & \cdots & \Delta s_{n+k} \\
\vdots & \vdots & \ddots & \vdots \\
\Delta^{k-1} s_n & \Delta^{k-1} s_{n+1} & \cdots & \Delta^{k-1} s_{n+k}
\end{vmatrix}
}{
\begin{vmatrix}
1 & 1 & \cdots & 1 \\
\Delta s_n & \Delta s_{n+1} & \cdots & \Delta s_{n+k} \\
\vdots & \vdots & \ddots & \vdots \\
\Delta^{k-1} s_n & \Delta^{k-1} s_{n+1} & \cdots & \Delta^{k-1} s_{n+k}
\end{vmatrix}
}.
\label{eq:shanks_determinant}
\end{equation}
The numerator contains the sequence values in the first row and forward differences in the subsequent rows while the denominator replaces the first row by ones.  
This determinant formulation is particularly elegant because it clearly shows the elimination of the first $k$ terms in an exponential error model.

\nd Assume the sequence has the exponential model
\begin{equation}
s_n = s^\ast + \sum_{j=1}^k c_j \lambda_j^n,
\end{equation}

\nd then the $k$-th order Shanks transformation exactly eliminates all $k$ terms, yielding
\[
e_k^{(n)} = s^\ast.
\]
For sequences with more than $k$ exponential components, $e_n^{(k)}$ eliminates the $k$ dominant terms, improving convergence to the next largest $\lambda_{k+1}$ component.

\nd The determinant form \eqref{eq:shanks_determinant} can also be interpreted in terms of Padé approximants. Let the generating function of $\{x_n\}$ be
\begin{equation}
f(z) = \sum_{n=0}^\infty c_n z^n.
\end{equation}

\nd The $k$-th order Shanks transform $e_n^{(k)}$ corresponds to evaluating the diagonal Padé approximant $[k/k]$ of $f(z)$ at $z=1$. This rational approximation viewpoint explains the method's effectiveness in accelerating sequences with multiple exponential terms or near singularities.

\begin{remark}
 The determinant denominators may become small, leading to round-off errors. Recursive $\varepsilon$-algorithm implementations often avoid this issue. Direct determinant evaluation is expensive for large $k$; recursion via $\varepsilon$-algorithm is more efficient.
 For sequences with unknown numbers of exponentials, iterated Shanks transformations provide robust convergence acceleration.
\end{remark}

\medskip 
\nd The general Shanks transformation is a powerful tool for convergence acceleration:
\begin{itemize}
\item The determinant form provides a **closed-form expression** for eliminating multiple exponential terms.
\item Recursive implementations ($\epsilon$-algorithm) make the method practical for numerical applications.
\item Connections with Padé approximants explain why the method often succeeds where linear extrapolation or Aitken $\Delta^2$ fails.
\item Applicable to sequences from numerical solutions of differential equations, iterative solvers, and series summation.
\end{itemize}

\medskip 
\nd Notice also that 
The Shanks transformation is mathematically equivalent to constructing diagonal Pad\'e approximants of the generating function $S(z)$. More precisely, the $k$-th Shanks transform corresponds to the Pad\'e approximant of type $(k,k)$.
\\
This equivalence explains why the Shanks transformation is particularly effective for sequences whose errors admit exponential representations of the form
\begin{equation}
s_n = s^\ast + \sum_{j=1}^k c_j \lambda_j^n,
\quad |\lambda_j| < 1.
\end{equation}
In such cases, Shanks transformation can eliminate all $k$ dominant exponential error terms exactly.\\
Shanks’ introduction of the transformation~(2.2) represented a major breakthrough
in the theory of sequence transformations. One of its most remarkable features is
the intimate connection with Pad\'e approximation, as demonstrated by Shanks
himself. 
\\
Let
\begin{equation}
f_n(z) = \sum_{\nu=0}^{n} \gamma_\nu z^\nu, \qquad n \in \mathbb{N}_0,
\label{eq:partialsums}
\end{equation}
denote the sequence of partial sums associated with the (formal) power series
expansion of a function $f$. According to a classical result (see~\cite{Shanks1955}),
the application of the Shanks transformation to these partial sums yields
\begin{equation}
e_k\!\left(f_n(z)\right) = [k+n\,/\,k]_f(z),
\qquad k,n \in \mathbb{N}_0,
\label{eq:shanks_pade}
\end{equation}
where $[k+n/k]_f(z)$ denotes the Pad\'e approximant of type $(k+n,k)$ associated
with the series of $f$.

\section{Wynn's $\varepsilon$-algorithm}

Wynn's $\varepsilon$-algorithm \cite{Wynn1956} is a powerful recursive method for accelerating the convergence of scalar sequences and summing slowly convergent or even divergent series. It generalizes the Shanks transformation to iteratively eliminate multiple exponential components in the sequence error and provides a practical, numerically stable implementation. 
\\
Consider a scalar sequence $\{x_n\}$ converging to a limit $s^\ast$:
\begin{equation}
s_n = s^\ast + \sum_{j=1}^\infty c_j \lambda_j^n, \quad |\lambda_1| > |\lambda_2| > \cdots.
\end{equation}

\nd The $\varepsilon$-algorithm systematically removes the leading exponential terms to accelerate convergence.
\\
Let $\varepsilon_{k}^{(n)}$ denote the $\varepsilon$-table elements. Initialize:
\begin{align}
\varepsilon_{-1}^{(n)} &= 0, \\
\varepsilon_{0}^{(n)} &= s_n.
\end{align}

\nd The recursion for $k \ge 0$ is:
\begin{equation}
\varepsilon_{k+1}^{(n)} = \varepsilon_{k-1}^{(n+1)} + \frac{1}{\varepsilon_{k}^{(n+1)} - \varepsilon_{k}^{(n)}}, \quad \varepsilon_{k}^{(n+1)} \neq \varepsilon_{k}^{(n)}.
\label{eq:epsilon_recursion}
\end{equation}

\nd The recursion proceeds along a triangular $\varepsilon$-table with $n$ indexing the sequence and $k$ the level of extrapolation.   Even-indexed elements $\varepsilon_{2k}^{(n)}$ are considered approximations to the sequence limit $x^\ast$, while odd-indexed elements are auxiliary and not directly used.
\\
Wynn's algorithm is equivalent to iteratively applying the Shanks transformation. Specifically:
\begin{equation}
\varepsilon_{2}^{(n)} = s_n^{(1)} = s_n - \frac{(\Delta s_n)^2}{\Delta^2 s_n},
\end{equation}
which is exactly the Aitken $\Delta^2$ process or first-order Shanks transformation.  

\nd Higher-order $\varepsilon_{2k}^{(n)}$ elements correspond to the $k$-th order Shanks transformation.   The algorithm provides a numerically stable and recursive implementation of the determinant formulas for the general Shanks transformation.

\nd The general $k$-th order Shanks transformation can be expressed in determinant form:
\begin{equation}
\varepsilon_{2k}^{(n)} = \frac{
\begin{vmatrix}
s_n & s_{n+1} & \cdots & s_{n+k} \\
\Delta s_n & \Delta s_{n+1} & \cdots & \Delta s_{n+k} \\
\vdots & \vdots & \ddots & \vdots \\
\Delta^{k-1} s_n & \Delta^{k-1} s_{n+1} & \cdots & \Delta^{k-1} s_{n+k}
\end{vmatrix}
}{
\begin{vmatrix}
1 & 1 & \cdots & 1 \\
\Delta s_n & \Delta s_{n+1} & \cdots & \Delta s_{n+k} \\
\vdots & \vdots & \ddots & \vdots \\
\Delta^{k-1} s_n & \Delta^{k-1} s_{n+1} & \cdots & \Delta^{k-1} s_{n+k}
\end{vmatrix}
}.
\end{equation}
The $\varepsilon$-algorithm evaluates these determinant forms recursively, avoiding direct determinant computations.  
This shows that $\varepsilon$-algorithm is the practical realization of the general Shanks transformation.
%\subsection{Convergence Properties}
\\
% Assume the sequence has the model:
% \begin{equation}
% x_n = x^\ast + \sum_{j=1}^m c_j \lambda_j^n, \quad |\lambda_1| > |\lambda_2| > \cdots > 0.
% \end{equation}
% Elimination of leading term: $\varepsilon_2^{(n)}$ removes the dominant term $c_1 \lambda_1^n$, yielding
% \[
% \varepsilon_2^{(n)} - x^\ast = \mathcal{O}(\lambda_2^n).
% \]
Notice that $\varepsilon_{2k}^{(n)}$ removes the first $k$ dominant exponential terms.   For sequences that are finite sums of exponentials with length $k$, the $\varepsilon_{2k}^{(n)}$ element is exactly equal to the limit $s^\ast$.

% \subsection{Application to Series Summation and Iterative Methods}

% 1. **Series acceleration:** For a slowly convergent series $S = \sum_{n=0}^\infty a_n$, define the partial sums $s_n = \sum_{i=0}^n a_i$ and apply the $\varepsilon$-algorithm to $\{s_n\}$.  

% 2. **Iterative methods:** For fixed-point iterations $x_{n+1} = g(x_n)$, apply the $\varepsilon$-algorithm to the sequence $\{x_n\}$ to accelerate convergence.  

% 3. **Connection with Padé approximants:** Applying $\varepsilon$-algorithm to the sequence of partial sums of a power series produces **diagonal Padé approximants** of the generating function.

% \subsection{Practical Considerations}

% - **Numerical stability:** The recursive formula avoids direct evaluation of large determinants. Care must be taken when $\varepsilon_k^{(n+1)} \approx \varepsilon_k^{(n)}$ to prevent division by nearly zero.  
% - **Storage:** Only two previous rows of the ε-table need to be stored at a time.  
% - **Stopping criteria:** Monitor convergence of even-indexed elements $\varepsilon_{2k}^{(n)}$.

\nd Wynn's $\varepsilon$-algorithm provides a powerful and numerically stable framework for sequence acceleration. It generalizes the Shanks transformation to higher order.  It has a recursive implementation that avoids determinant evaluation.  It can accelerate sequences, partial sums of series, and iterative fixed-point methods.   It has also  a natural connection to Padé approximants, explaining its effectiveness in accelerating convergence of sequences with exponential or near-exponential error components. For more details on convergence of the $\varepsilon$-algorithm, see \cite{Brezinski1971r,Wynn1956,Wynn1962,Wynn1966}.

% \section{Levin-Type Transformations}

% Levin-type methods exploit explicit remainder estimates:
% \[
% x_n = x^\ast + \omega_n \sum_{k=0}^{\infty} c_k \phi_k(n),
% \]
% where $\omega_n$ is a known remainder estimate.

% The Levin $u$-transformation is given by
% \begin{equation}
% u_n =
% \frac{\sum_{k=0}^m (-1)^k \binom{m}{k}
% (n+k)^{m-1} x_{n+k} / \omega_{n+k}}
% {\sum_{k=0}^m (-1)^k \binom{m}{k}
% (n+k)^{m-1} / \omega_{n+k}}.
% \end{equation}

% These methods are highly effective for slowly convergent or divergent series.

% \section{Numerical Stability and Conditioning}

% Despite their efficiency, scalar extrapolation methods can suffer from:
% \begin{itemize}
% \item Cancellation errors,
% \item Division by small differences,
% \item Sensitivity to perturbations and noise.
% \end{itemize}

% Stabilization strategies include adaptive selection of extrapolation depth and hybrid methods.

% \section{Applications}

% Scalar extrapolation techniques are widely used in:
% \begin{itemize}
% \item Fixed-point iterations,
% \item Numerical quadrature,
% \item Evaluation of special functions,
% \item Summation of divergent series,
% \item Error estimation and adaptivity.
% \end{itemize}

% \section{Recent Developments}

% \nd Recent research directions include:
% \begin{itemize}
% \item Robust extrapolation under stochastic noise,
% \item Adaptive and automatic extrapolation schemes,
% \item Combination with vector and tensor acceleration methods,
% \item Applications in scientific computing and data science.
% \end{itemize}

\subsection{Connection with Pad\'e approximation}

\nd One of the most profound theoretical interpretations of Wynn's $\varepsilon$-algorithm is its intimate connection with Pad\'e approximation. This relationship provides a unifying framework linking sequence transformations, rational approximation, and convergence acceleration.\\
\emph{Pad\'e Approximation}. Consider the function corresponding to the sum of the formal power series
\begin{equation}
f(z) = \sum_{k=0}^{\infty} a_k z^k
\label{eq:power-series}
\end{equation}
The $(m,n)$ Pad\'e approximant $[m/n]_f(z)$ is the rational function
\begin{equation}
[m/n]_f(z)
=
\frac{P_m(z)}{Q_n(z)},
\qquad
Q_n(0) = 1,
\label{eq:pade-def}
\end{equation}
where $P_m$ and $Q_n$ are polynomials of degrees at most $m$ and $n$, respectively, chosen such that
\begin{equation}
f(z) - [m/n]_f(z)
=
\mathcal{O}(z^{m+n+1})
\quad \text{as } z \to 0.
\label{eq:pade-cond}
\end{equation}
The collection of Pad\'e approximants $\{[m/n]\}$ arranged in a two-dimensional array indexed by $(m,n)$ is known as the \emph{Pad\'e table}; for more details on Padé approximation, refer to \cite{Baker1981,Brezinski1980}.\\
Let $(s_n(z))_n$ be the sequence defined by 
$$s_n(z)=\displaystyle \sum_{k=0}^n a_k z^k.$$
Applying Wynn’s $\varepsilon$-algorithm to the sequence $\{s_n(z)\}$, viewed as a function of the parameter $z$, yields the fundamental identity
\begin{equation}
\varepsilon_{2k}^{(n)}(z)
=
[n+k/k]_{f}(z),
\label{eq:eps-pade}
\end{equation}
provided the relevant Pad\'e approximant exists and the algorithm does not encounter a breakdown.
\\
Thus, the even elements of the $\varepsilon$-table coincide exactly with the diagonal and near-diagonal entries of the Pad\'e table. In particular,
\[
\varepsilon_{2k}^{(0)}(z) = [k/k]_f(z),
\]
corresponding to the main diagonal of the Pad\'e table. Thus, the $\varepsilon$-algorithm may be viewed as a \emph{Pad\'e table generator}, with the index $k$ corresponding to the order of the Pad\'e approximant and the index $n$ acting as a shift along the sequence. 
This Pad\'e interpretation explains the remarkable efficiency of the $\varepsilon$-algorithm. While polynomial approximations reproduce only algebraic behavior, rational approximations are capable of representing functions with poles, branch points, or other singularities. Consequently, diagonal Pad\'e approximants often converge in regions where the original power series diverges or converges very slowly.

\nd From the viewpoint of sequence transformations, the Pad\'e connection implies that the $\varepsilon$-algorithm eliminates dominant exponential error terms in the sequence model
\begin{equation}
s_n = s + \sum_{j=1}^{m} c_j \lambda_j^n,
\end{equation}
by constructing a rational approximant whose denominator captures the dominant poles associated with the $\lambda_j$. As a result, each increase in $k$ removes an additional exponential component, leading to dramatically accelerated convergence.
%\paragraph{Relation to Shanks Transformation.}

\nd The equivalence between Pad\'e approximants and the Shanks transformation becomes explicit through determinant representations. The general Shanks transformation produces exactly the same approximations as diagonal Pad\'e approximants, while the $\varepsilon$-algorithm provides a numerically stable recursive implementation of these determinant formulas. In this sense, the $\varepsilon$-algorithm can be interpreted as an efficient computational mechanism for evaluating Pad\'e approximants without explicitly constructing numerator and denominator polynomials.
\\
%\paragraph{Practical Implications.}
The Pad\'e-table interpretation has several important consequences: It explains why the $\varepsilon$-algorithm is particularly effective for summing slowly convergent or divergent series,  it provides a theoretical justification for its robustness in the presence of singularities in the generating function and it establishes a direct link between convergence acceleration techniques and classical rational approximation theory.

\nd Viewed as a Pad\'e table, Wynn's $\varepsilon$-algorithm occupies a central position at the intersection of sequence transformations, rational approximation, and numerical analysis. This interpretation not only clarifies its convergence properties but also explains its enduring success in practical computations involving series summation and iterative processes.

 \subsection{Continued fractions and rational approximation}

 Continued fractions provide a natural and powerful framework for understanding rational approximation and play a central role in the theoretical analysis of extrapolation and convergence acceleration methods. In particular, Pad\'e approximants admit representations in terms of continued fractions, and many extrapolation algorithms can be interpreted as procedures that implicitly generate such representations.\\

 \nd\emph{Pad\'e Approximants and Continued Fractions.}
 Given a formal power series
\begin{equation}
f(z) = \sum_{n=0}^{\infty} c_n z^n,
\end{equation}
and consider the continued fraction
\begin{equation}
 C(z)=b_0+\cfrac{a_1z}{b_1 + \cfrac{a_2 z}{b_2 + \cfrac{a_3 z}{b_3 + \ddots}}}.
\end{equation}
Let $C_n(z)=
 \displaystyle \frac{A_n(z)}{B_n(z)}$ denote the $n$-th convergent  of $C(z)$. 
The continued fraction converges for a fixed $z$ if and only if $(C_n(z))_n$ converges to  $C(z)$ as $n$ goes to infinity.  Pad\'e approximants of $f$ may be expressed as truncations of a continued fraction and we have 
$$C_{2k}=[k/k]_f(z)\;\; \text{and}\;\;C_{2k+1}=[k+1/k]_f(z). $$
Notice also that to obtaining the partial numerator and denominator  of the corresponding continued fraction, one can use the well knwon \emph{qd-algorithm} of H. Rutishauser \cite{Rutishauser1957}.

\nd \emph{Implicit Generation by the $\varepsilon$-Algorithm.}
Although Wynn's $\varepsilon$-algorithm does not explicitly construct continued fractions, it implicitly generates the coefficients of such a representation through its recursive updates. Each even column of the $\varepsilon$-table corresponds to a successive convergent of a continued fraction associated with the generating function of the sequence. In this sense, the $\varepsilon$-algorithm may be interpreted as a \emph{continued fraction accelerator}, transforming partial sums into rational approximants that capture the dominant analytic structure of the underlying function.\\

\subsection{The $E$-algorithm}

Brezinski’s $E$-algorithm \cite{BrezinskiR1991}, is a general sequence transformation designed to accelerate convergence by eliminating prescribed components of the asymptotic error expansion. It provides a unifying framework that encompasses a wide class of extrapolation methods, including Richardson extrapolation, the Shanks transformation, Levin-type transformations, and Wynn’s $e$-algorithm. Unlike purely recursive schemes such as Wynn’s $\varepsilon$-algorithm, the $E$-algorithm is rooted in an explicit algebraic formulation based on annihilation of basis functions modeling the error.
\\
Let $\{s_n\}$ be a scalar sequence converging to a limit $s$, and assume that the truncation error admits an expansion of the form
\begin{equation}
s_n = s + \sum_{j=1}^{\infty} c_j g_j(n),
\end{equation}
where $\{g_j(n)\}$ is a sequence of known functions.
%such that
% \[
% \lim_{n\to\infty} g_{j+1}(n)/g_j(n) = 0.
% \]
The functions $g_j(n)$ describe the dominant asymptotic behavior of the error and may represent algebraic, exponential, or logarithmic decay.\\
%\paragraph{Definition of the $E$-Algorithm.}
For a fixed integer $k$, the $E$-algorithm constructs an approximation $E_k^{(n)}$ to the limit $s$ by imposing the exactness condition
\begin{equation}
E_k^{(n)} = s \quad \text{whenever} \quad
s_n = s + \sum_{j=1}^{k} c_j g_j(n).
\end{equation}
This leads to a linear system of equations whose solution yields the transformed new approximation 
\begin{equation}
E_k^{(n)} =
\frac{
\begin{vmatrix}
s_n & s_{n+1} & \cdots & s_{n+k} \\
g_1(n) & g_1(n+1) & \cdots & g_1(n+k) \\
\vdots & \vdots &  & \vdots \\
g_k(n) & g_k(n+1) & \cdots & g_k(n+k)
\end{vmatrix}
}{
\begin{vmatrix}
1 & 1 & \cdots & 1 \\
g_1(n) & g_1(n+1) & \cdots & g_1(n+k) \\
\vdots & \vdots &  & \vdots \\
g_k(n) & g_k(n+1) & \cdots & g_k(n+k)
\end{vmatrix}
}.
\end{equation}
This determinant expression generalizes the Shanks transformation and makes explicit the algebraic structure underlying the method.\\

\begin{theorem}
    $\forall n$, $E_k^{(n)}= S$ if and only if
    $$s_n=s+a_1g_1(n)+\ldots+ a_kg_k(n),\; k, n \ge 0.$$
\end{theorem}
Let us write now $E_k^{(n)}$ as a Schur complement. We denote by ${\bf g}_n$ the vector  ${\bf g}_n=(g_1(n),\ldots,g_k(n))^T$, ${\bf G}_k^{(n)}$ the $k \times k$ 
$$[{\bf G}_k^{(n)}]_{ij}=\Delta g_i(n+j-1),\; i,j=1,\ldots,k. $$
% matrix whose columns are the vectors $(\Delta g_1(n),\ldots,\Delta g_k(n))^T,\ldots, (\Delta g_1(n+k-1),\ldots,\Delta g_k(n+k-1))^T$. 
Then, $E_k^{((n)}$ can also be expressed as
$$E_k^{((n)}=
\begin{vmatrix}
s_n& \Delta S_{k,n}\\
{\bf g}_n & {\bf G}_k^{(n)}
    \end{vmatrix} / \begin{vmatrix}
        {\bf G}_k^{(n)},
    \end{vmatrix}$$
    where $\Delta S_{k,n}=[\Delta s_n,\ldots, \Delta s_{n+k-1} ].$
Then using the Schur complement formulation and assuming that the matrix ${\bf G}_k^{(n)}$ is nonsingular, we obtain:
$$E_k^{(n)}= s_n -\Delta S_{k,n}\, {\left({\bf G}_k^{(n)}\right)}^{-1} {\bf g}_n.$$

\nd \emph{Recursive Implementation.}
Although defined via determinants, the $E$-algorithm admits stable recursive implementations that avoid direct determinant evaluation. These recursions significantly reduce computational cost and improve numerical robustness, making the method suitable for practical applications.
\\
The $E$‑algorithm, introduced by Brezinski, provides a general framework for accelerating the convergence of sequences by explicitly modeling and annihilating dominant error components. While the theoretical definition uses determinants or exact annihilation conditions, a practical and efficient implementation relies on a recursive scheme that avoids computationally expensive determinants and improves numerical stability.
\\

\nd The core of the recursive implementation lies in constructing a two‑indexed table $(\mathbf{E}_k^{(n)})$ defined as follows:
$$
\left\{
\begin{array}{lll}
\mathbf{E}_{0}^{(n)} = \mathbf{s}_n,\;  g_{0,i}(n)=g_i(n),\\
\mathbf{E}_{k}^{(n)} = 
\mathbf{E}_{k-1}^{(n+1)}
 - \displaystyle \frac{\mathbf{E}_{k-1}^{(n+1)}-\mathbf{E}_{k-1}^{(n)}}{\mathbf{g}_{k-1,k}^{(n+1)}-\mathbf{g}_{k-1,k}^{(n)}} {\mathbf{g}_{k-1,k}^{(n)}}\\
 \mathbf{g}_{k,i}^{(n)} = 
\mathbf{g}_{k-1,i}^{(n+1)}
 - \displaystyle \frac{\mathbf{g}_{k-1,i}^{(n+1)}-\mathbf{g}_{k-1,i}^{(n)}}{\mathbf{g}_{k-1,k}^{(n+1)}-\mathbf{g}_{k-1,k}^{(n)}} \mathbf{g}_{k-1,k}^{(n)}, 
\quad i=k+1,\ldots, 
\end{array}
\right.
$$ 
with the convention 
$ \mathbf{E}_{-1}^{(n)} = \mathbf{0},\quad \forall n$.
\\
The sequence $g_{k,i}^{(n)}$ has also the following determinant  expression
\begin{equation}
g_{k,i}^{(n)} =
\frac{
\begin{vmatrix}
g_i(n) & g_i(n+1) & \cdots & g_i(n+k) \\
g_1(n) & g_1(n+1) & \cdots & g_1(n+k) \\
\vdots & \vdots &  & \vdots \\
g_k(n) & g_k(n+1) & \cdots & g_k(n+k)
\end{vmatrix}
}{
\begin{vmatrix}
1 & 1 & \cdots & 1 \\
g_1(n) & g_1(n+1) & \cdots & g_1(n+1) \\
\vdots & \vdots &  & \vdots \\
g_k(n) & g_k(n+1) & \cdots & g_k(n+k)
\end{vmatrix}
}.
\end{equation}

\medskip 
Here also, $g_{k,i}^{(n)}$ can be expressed as the following Schur complement
$$g_{k,i}^{(n)}= g_i(n) -[\Delta g_i(n),\ldots, \Delta g_i(n+k-1) ]\, {\left({\bf G}_k^{(n)}\right)}^{-1} {\bf g}_n.$$

\nd Different choices of the basis functions $g_j(n)$ yield well-known extrapolation techniques such as  Shanks transformation,  Richardson extrapolation and polynomial or rational choices give rise to generalized Pad\'e-type approximants. 
Thus, the $E$-algorithm serves as a unifying algebraic framework for convergence acceleration, for more details, refer to  the book \cite{BrezinskiR1991}. The  $E$-algorithm is a powerful and general framework for sequence transformation based on explicit asymptotic modeling. By unifying a wide range of classical and modern extrapolation methods, it plays a central role in the theory of convergence acceleration and provides deep insight into the algebraic and analytic structure of extrapolation techniques.
Other well know scalar extrapolation methods are developed in \cite{BrezinskiR1991}. Among them the $\rho$ and $\theta$-algorithms that we will describe in the following.

\subsection{The $\rho$ algorithm}

The iterated $\Delta^2$ process  and Wynn’s $\varepsilon$-algorithm 
belong to the most effective classical tools for accelerating the convergence of
linearly convergent sequences. Beyond this, they are also capable of assigning finite
values to a wide class of alternating divergent series. 
% Comprehensive discussions and
% applications can be found, for instance, in~\cite[Section~15.2]{BrezinskiBook}, as well
% as in the rigorous convergence analysis of the Euler series presented in~\cite{HardyEuler}.
\\
Despite their broad applicability, both the iterated $\Delta^2$ process and the
$\varepsilon$-algorithm fail when applied to sequences exhibiting logarithmic
convergence. This limitation is fundamental and has been formally established by some authors. In order to overcome this deficiency, Wynn introduced in \cite{Wynn1956} a new extrapolation method,  known as the \emph{$\rho$-algorithm} .
\\
The $\rho$-algorithm is defined recursively by
$$
\left\{ 
\begin{array}{ll}
\rho_{-1}^{(n)} = 0, \qquad
\rho_{0}^{(n)} = s_n, \qquad n \in \mathbb{N},
\\
\rho_{k+1}^{(n)} =
\rho_{k-1}^{(n+1)} +
\displaystyle \frac{x_{n+k+1} - x_n}
{\rho_k^{(n+1)} - \rho_k^{(n)}},
\qquad k,n \in \mathbb{N}.
\end{array}
\right.,
$$
where $x_0,x_1,\ldots$ are given numbers. 
This algorithm is particularly
well suited for accelerating the convergence of sequences whose asymptotic behavior
is governed by logarithmic terms. 
Similarly to the $\varepsilon$-algorithm, the $\rho$-algorithm produces meaningful
approximations to the limit of the original sequence only through its elements with
even subscripts. That is, the quantities $\rho_{2k}^{(n)}$ represent the extrapolated
approximations to the limit of the sequence $\{s_n\}$, whereas the elements $\rho_{2k+1}^{(n)}$
serve merely as auxiliary quantities. If the transformation converges, these odd-index
elements typically diverge and should not be interpreted as approximations.
\\
From a theoretical perspective, the quantities $\rho_{2k}^{(n)}$ admit a natural
interpretation in terms of continued fractions. More precisely, they coincide with
the terminants of a continued fraction that interpolates the sequence $\{s_n\}$ at a
set of nodes $\{x_n\}$ and extrapolates this interpolation process to infinity. This interpretation immediately
imposes structural conditions on the interpolation points $\{x_n\}$, which must form
a strictly increasing and unbounded sequence:
\begin{equation}
0 < x_0 < x_1 < \cdots < x_m < x_{m+1} < \cdots,
\label{eq:xn_monotone}
\end{equation}
and
\begin{equation}
\lim_{n \to \infty} x_n = \infty.
\label{eq:xn_infty}
\end{equation}
In practical computations, Wynn’s $\rho$-algorithm is most often employed in its
standard form, obtained by choosing the interpolation points as
\[
x_n = n+1, \qquad n \in \mathbb{N}.
\]
With this choice, the recursion simplifies to
$$
\left\{
\begin{array}{ll}
\rho_{-1}^{(n)} = 0, \qquad
\rho_{0}^{(n)} = s_n, \\
\rho_{k+1}^{(n)} =
\rho_{k-1}^{(n+1)} +
\displaystyle \frac{k+1}
{\rho_k^{(n+1)} - \rho_k^{(n)}},
\qquad k,n \in \mathbb{N}.
\label{eq:rho_standard}
\end{array}
\right.
$$
This formulation is the one most frequently encountered in the literature and in
numerical applications; see for example ~\cite{Brezinski1975}.

\subsubsection{The $\theta$-algorithm}

The $\theta$-algorithm was introduced by Brezinski in \cite{Brezinski1971} as a nonlinear sequence
transformation designed to overcome certain limitations of Wynn’s
$\varepsilon$-algorithm and $\rho$-algorithm. In particular, while the
$\varepsilon$-algorithm is highly effective for linearly convergent sequences
and the $\rho$-algorithm is well adapted to logarithmically convergent ones,
neither method performs satisfactorily when the sequence exhibits a mixture of
linear and logarithmic convergence behaviors. The $\theta$-algorithm was
constructed precisely to address this intermediate situation.
\\
Let $\{s_n\}_{n \ge 0}$ be a scalar sequence. The $\theta$-algorithm is defined
by the recursive scheme
$$
\left\{ 
\begin{array}{lll}
\theta_{-1}^{(n)} = 0, \qquad
\theta_{0}^{(n)} = s_n, \qquad n \in \mathbb{N},
\\
\theta_{2k+1}^{(n)} =
\theta_{2k-1}^{(n+1)} + D_{2k}^{(n)}\\
\theta_{2k+2}^{(n)} =  \theta_{2k}^{(n+1)} - \displaystyle \frac{\Delta \theta_{2k}^{(n+1)}}{\Delta D_{2k+1}^{(n)}}  D_{2k+1}^{(n)}
\qquad k,n \in \mathbb{N},
\end{array}
\right.
$$
provided that the denominator does not vanish.
\\
As in the case of the $\varepsilon$- and $\rho$-algorithms, only the elements
$\theta_{2k}^{(n)}$ with even subscripts are intended to approximate the limit
of the original sequence. The quantities $\theta_{2k+1}^{(n)}$ serve as auxiliary
terms and typically diverge when the transformation converges.

\medskip
\noindent
The defining feature of the $\theta$-algorithm is the explicit use of second
forward differences in the denominator. This
structure allows the transformation to detect and exploit curvature in the
sequence, which is characteristic of logarithmic convergence, while retaining
the ability to accelerate linear convergence. As a consequence, the
$\theta$-algorithm often succeeds in cases where both the $\varepsilon$- and
$\rho$-algorithms fail or stagnate.
\\
From a theoretical standpoint, the $\theta$-algorithm may be viewed as a hybrid
method that combines ideas from Aitken-type acceleration and Wynn-type nonlinear
recursions. It does not rely on external interpolation points, in contrast to
the $\rho$-algorithm, and it does not require prior knowledge of the asymptotic
form of the remainder.

\medskip
\noindent
If the sequence $\{s_n\}$ converges linearly, the $\theta$-algorithm exhibits
convergence properties comparable to those of the $\varepsilon$-algorithm. For
logarithmically convergent sequences, it often provides a significant
improvement over both the $\varepsilon$- and $\rho$-algorithms. This robustness
has made the $\theta$-algorithm a valuable tool in numerical analysis,
particularly in the summation of slowly convergent series and in extrapolation
procedures arising in applied mathematics and physics. 
Further theoretical analysis, extensions, and numerical examples of the
$\theta$-algorithm can be found in the works of Brezinski and collaborators.

\section{Vector extrapolation methods}
The development and success of sequence transformation and extrapolation techniques  allow the introduction of a variety of vector sequence extrapolation
methods since 1970s to address problems arising in numerical linear algebra and fixed-point iterations. Notable among these vector extrapolation methods are the minimal
polynomial extrapolation (MPE) method of Cabay and Jackson \cite{CabyJ1976}, the
reduced rank extrapolation (RRE) method proposed by Eddy \cite{Eddy1979} and
Mesina \cite{Messina1977}, and the modified minimal polynomial extrapolation
(MMPE) method developed by Pugatchev \cite{Pugatchev1978} and further analyzed by
Brezinski \cite{Brezinski1975} and the vector $\varepsilon$-algorithm of Wynn \cite{Wynn1962}. Comprehensive theoretical treatments and practical
implementations of these methods can be found in 
\cite{Brezinski1991,Jbilousadok2000,Salam1996,Sidi2003,Wynn1962}.
\\
The polynomial extrapolation framework was further clarified and extended through
connections with Padé approximation and orthogonal polynomials
\cite{Brezinski1980}, as well as through detailed algorithmic studies for linear
systems and fixed-point problems. Significant contributions in this direction
were made in \cite{BrezinskiR2014,BrezinskiCRS2022,Jbilousadok1995,Jbilousadok1999,Jbilousadokreichel2009,Salam1996,Sidi2003}. Related advances include the works
in  \cite{Sidi1991,SidiFS1988} and applications to rank-deficient and ill-conditioned problems discussed by Hansen \cite{Hansen1994,Hansen1998}.
\\
More recently, these ideas have been generalized beyond vectors to block and
tensor settings. Block extrapolation methods were investigated in \cite{Jbiloumessaoudi2016}, while tensor extrapolation techniques and their applications were developed in \cite{BeikElJS2021,BentbibJT2024}. These extensions demonstrate that extrapolation
methods remain a flexible and powerful tool for accelerating convergence in
high-dimensional and structured numerical problems. Vector extrapolation methods could be divided into two categories: the polynomial   ones including MPE, RRE and MMPE and $\varepsilon$-based ones containing vector and topological $\varepsilon$-algorithms. \\
Since the Schur complement formula will be used extensively in the sequel, we
recall its definition here for completeness.\\Let
\[
M =
\begin{pmatrix}
A & B \\
C & D
\end{pmatrix},
\]
where $A \in \mathbb{R}^{m \times m}$ and $D \in \mathbb{R}^{n \times n}$ are
square matrices. If $A$ is invertible, the \emph{Schur complement} of $A$ in $M$
is defined by
\[
(M/A) = D - C A^{-1} B.
\]
Similarly, if $D$ is invertible, the Schur complement of $D$ in $M$ is given by
\[
(M/D) = A - B D^{-1} C.
\]
Notice for example  that $det[(M/A)]=det(M)/det(A).$

\subsection{The polynomial-based methods}
\subsubsection{Vector Polynomial Extrapolation (VPE) methods}

In this subsection, we describe some of the most widely used vector polynomial
extrapolation (VPE) methods, namely the Minimal Polynomial Extrapolation (MPE) method
\cite{CabyJ1976}, the Reduced Rank Extrapolation (RRE) method
\cite{Eddy1979,Messina1977}, and the Modified Minimal Polynomial Extrapolation
(MMPE) method \cite{Brezinski1975,Pugatchev1978}. These methods are designed to
accelerate the convergence of a sequence of vectors $(s_n) \subset \mathbb{R}^N$
towards a limit (or a fixed point) by forming suitable linear combinations of
the sequence terms.
\\
Let $(s_n)$ be a given sequence of vectors, and consider a general transformation
$T_k$ of order $k$ defined as
\[
T_k: \mathbb{R}^N \longrightarrow \mathbb{R}^N, \quad s_n \mapsto t_k^{(n)},
\]
with
\begin{equation}
t_k^{(n)} = s_n + \sum_{i=1}^k a_i^{(n)} \phi_i(n), \quad n \ge 0,
\end{equation}
where $(\phi_i(n))_n$ are auxiliary vector sequences chosen according to the particular
extrapolation method. The idea is to construct $t_k^{(n)}$ as a correction to
$s_n$, based on a linear combination of increments or differences of the
sequence.
\\
We also define the shifted transformation $\widetilde T_k$ by
\begin{equation}
\tilde t_k^{(n)} = s_{n+1} + \sum_{i=1}^k a_i^{(n)} \phi_i(n+1), \quad n \ge 0,
\end{equation}
where the coefficients $a_i^{(n)}$ are identical to those used in $T_k$. This
shifted sequence allows us to define the \emph{generalized residual} of
$t_k^{(n)}$:
\begin{align}
\tilde r(t_k^{(n)}) &= \tilde t_k^{(n)} - t_k^{(n)} \nonumber \\
&= \Delta s_n + \sum_{i=1}^k a_i^{(n)} \Delta \phi_i(n),
\end{align}
where the forward difference operator $\Delta$ acts on the sequence index:
\[
\Delta \phi_i(n) = \phi_i(n+1) - \phi_i(n), \quad i=1,\ldots,k.
\]
For polynomial vector extrapolation methods, the auxiliary sequences are usually chosen as
\begin{equation}
\phi_i(n) = \Delta s_{n+i-1}, \quad i=1,\ldots,k,
\end{equation}
which means that the correction terms are based on successive differences of the
sequence. The coefficients $a_i^{(n)}$ are then determined by imposing an
orthogonality condition on the generalized residual:
\begin{equation}
\tilde r(t_k^{(n)}) \perp \mathrm{span}\{y_1^{(n)}, \ldots, y_k^{(n)}\},
\end{equation}
where the choice of $y_i^{(n)}$, $i=1,\ldots,k$,  distinguishes the methods:
\[
y_i^{(n)} =
\begin{cases}
\Delta s_{n+i-1},\; \text{for}\;  & \text{MPE} \\
\Delta^2 s_{n+i-1} ,\; \text{for}\; & \text{RRE} \\
y_i\; \text{fixed vectors},\; \text{for}\; & \text{MMPE}.
\end{cases}
\]
Here, $\Delta^2 s_n = \Delta (\Delta s_n)$ denotes the second forward difference,
and $y_i$ in MMPE is typically chosen according to a minimal polynomial
construction.
\\
Let us define the subspaces
\[
\widetilde W_{k,n} = \mathrm{span}\{\Delta^2 s_n, \ldots, \Delta^2 s_{n+k-1}\}, \quad
\widetilde L_{k,n} = \mathrm{span}\{y_1^{(n)}, \ldots, y_k^{(n)}\}.
\]
Then, the generalized residual satisfies
\begin{align}
\tilde r(t_k^{(n)}) - \Delta s_n &\in \widetilde W_{k,n},\\
\tilde r(t_k^{(n)}) &\perp \widetilde L_{k,n}.
\end{align}
In matrix form, let $L_{k,n}, \Delta S_{k,n}, and \Delta^2 S_{k,n}$ be the $N \times k$
matrices whose columns are $y_1^{(n)}, \dots, y_k^{(n)}$, $\Delta s_n, \dots,
\Delta s_{n+k-1}$, and $\Delta^2 s_n, \dots, \Delta^2 s_{n+k-1}$, respectively.
Then
\begin{equation}
\tilde r(t_k^{(n)}) = \Delta s_n - \Delta^2 S_{k,n} \big(L_{k,n}^T \Delta^2 S_{k,n}\big)^{-1} L_{k,n}^T \Delta s_n,
\end{equation}
and the approximation $t_k^{(n)}$ is given by
\begin{equation}
t_k^{(n)} = s_n - \Delta S_{k,n} \big(L_{k,n}^T \Delta^2 S_{k,n}\big)^{-1} L_{k,n}^T \Delta s_n,
\end{equation}
which can also be expressed as a Schur complement. In fact, let ${\cal T}_k^{(n)}$ be the matrix defined by
\begin{equation}\label{sch1}
{\cal T}_k^{(n)}= 
\begin{pmatrix}
s_n & \Delta S_{k,n}\\
&\\
L_{k,n}^T \Delta s_n& L_{k,n}^T \Delta^2 S_{k,n}
\end{pmatrix}.
\end{equation}
Hence the approximation $t_k^{(n)}$ can be given as the following Schur complement
\begin{equation}\label{sch2}
t_k^{(n)} = \left( {\cal T}_k^{(n)} / L_{k,n}^T \Delta^2 S_{k,n} \right).
\end{equation}
\nd  Notice that $t_k^{(n)}$ can also be written as the following ratio of two determinants
\begin{equation}\label{sch3} 
 t_k^{(n)}=  
 \begin{vmatrix}
 s_n & \Delta S_{k,n}\\
 &\\
 L_{k,n}^T \Delta s_n & L_{k,n}^T \Delta^2 S_{k,n}
\end{vmatrix}/ \begin{vmatrix}
    L_{k,n}^T \Delta^2 S_{k,n},
\end{vmatrix}
\end{equation}
where the vector determinant in the numerator is the vector obtained by expanding it with respect to its first row.\\
Equivalently, $t_k^{(n)}$ can also be expressed as a convex linear combination of
sequence elements (see \cite{Jbilousadok2000}):
\begin{equation}
t_k^{(n)} = \sum_{j=0}^k \beta_j^{(n)} s_{n+j}, \quad \sum_{j=0}^k \beta_j^{(n)} = 1,
\end{equation}
where the coefficients satisfy the system
\begin{equation}
\sum_{j=0}^k \alpha_{i,j}^{(n)} \beta_j^{(n)} = 0, \quad i=0,\ldots,k-1,
\end{equation}
with
\[
\alpha_{i,j}^{(n)} =
\begin{cases}
(\Delta s_{n+i}, \Delta s_{n+j}), & \text{MPE}, \\
(\Delta^2 s_{n+i}, \Delta s_{n+j}), & \text{RRE}, \\
(y_{i+1}, \Delta s_{n+j}), & \text{MMPE}.
\end{cases}
\]
Finally, $t_k^{(n)}$ admits also  the following  determinant representation:
\begin{equation}
t_k^{(n)} =
\frac{
\begin{vmatrix}
s_n & s_{n+1} & \cdots & s_{n+k} \\
\alpha_{0,0}^{(n)} & \alpha_{0,1}^{(n)} & \cdots & \alpha_{0,k}^{(n)} \\
\vdots & \vdots & & \vdots \\
\alpha_{k-1,0}^{(n)} & \alpha_{k-1,1}^{(n)} & \cdots & \alpha_{k-1,k}^{(n)}
\end{vmatrix}
}{
\begin{vmatrix}
1 & 1 & \cdots & 1 \\
\alpha_{0,0}^{(n)} & \alpha_{0,1}^{(n)} & \cdots & \alpha_{0,k}^{(n)} \\
\vdots & \vdots & & \vdots \\
\alpha_{k-1,0}^{(n)} & \alpha_{k-1,1}^{(n)} & \cdots & \alpha_{k-1,k}^{(n)}
\end{vmatrix}
}.
\end{equation}
This representation highlights the connection between vector extrapolation methods
and classical determinants, making explicit the linear dependencies among the
sequence differences that are exploited to accelerate convergence.\\
It is clear that RRE is an orthogonal projection while MPE and MMPE are oblique projection methods. Reduced Rank Extrapolation constructs the extrapolated vector by minimizing the norm of a linear combination of forward differences. Specifically, RRE solves the constrained least-squares problem
\begin{equation}
\min_{\boldsymbol{\alpha}} 
\left\|
\sum_{j=0}^{k} \alpha_j \Delta s_{n+j}
\right\|,
\quad \text{subject to} \quad \sum_{j=0}^{k} \alpha_j = 1.
\end{equation}
The extrapolated vector is then given by
\begin{equation}
\mathbf{x}_n^{\text{RRE}} = \sum_{j=0}^{k} \alpha_j  s_{n+j}.
\end{equation}
RRE is numerically robust and equivalent to the GMRES method when applied to linear fixed-point iterations. Its least-squares formulation ensures stability and makes it particularly suitable for large-scale problems.

\medskip
\nd Minimal Polynomial Extrapolation determines the coefficients $\{\alpha_j\}$ by enforcing exact annihilation of a linear combination of the forward differences:
\begin{equation}
\sum_{j=0}^{k} \alpha_j \Delta s_{n+j} = \mathbf{0},
\quad \sum_{j=0}^{k} \alpha_j = 1.
\end{equation}
The resulting extrapolated vector is
\begin{equation}
\mathbf{x}_n^{\text{MPE}} = \sum_{j=0}^{k} \alpha_j s_{n+j}.
\end{equation}
The MPE method yields an exact result whenever the error can be represented by an exponential expansion involving no more than \(k\) components. Nevertheless, since its implementation requires the solution of a linear system that may be ill-conditioned, MPE can exhibit reduced numerical robustness compared with RRE when computations are carried out in finite-precision arithmetic. When applied to nonlinear problems, the behavior of MPE differs noticeably from that of RRE and Anderson acceleration. While MPE can still be effective when the
nonlinear iteration behaves approximately linearly near the solution, its performance is generally more sensitive to deviations from linearity

\nd MMPE imposes the orthogonality  condition on a reduced subspace spanned by selected test vectors $\{{y}_i\}_{i=1}^{k}$:
\begin{equation}
\sum_{j=0}^{k} \alpha_j\,  y_i^T \Delta s_{n+j}  = 0,
\quad i = 1, \ldots, k.
\end{equation}
This leads to a smaller linear system for the coefficients $\{\alpha_j\}$. The extrapolated vector is again given by
\begin{equation}
\mathbf{x}_n^{\text{MMPE}} = \sum_{j=0}^{k} \alpha_j s_{n+j}.
\end{equation}
In practical computations, the selection among RRE, MPE, and MMPE is guided by the balance between numerical stability and
exactness. RRE is often favored in large-scale settings because of its robustness, whereas MPE can be more effective when the
problem is well conditioned and high accuracy is sought. MMPE is an intermediate position between these two approaches.
Overall, polynomial vector extrapolation methods form a flexible and powerful framework for accelerating convergence. By
leveraging polynomial error representations and Krylov subspace properties, methods such as RRE, MPE, and MMPE have proved
effective in a broad range of scientific and engineering applications 
\cite{BrezinskiCRS2022,DuminilSS,GanderGG}.

\subsubsection{The $S\beta$-algorithm}
Let $\{s_n\}_{n\ge0}$ be a sequence of vectors in $\mathbb{R}^N$. Let $k$ be a positive integer such that $k<p$, and let
$n\in\mathbb{N}$. We consider the transformation $(s_{k}^n)$ obtained by the MMPE method with the given independent vectors $y_1,\ldots,y_k$ and let $Y_k$ the matrix whose columns are those vectors: 
\begin{equation}
s_{k}^{(n)}
=
s_n
-
\Delta S_{n,k}
\bigl(
Y_{k}^T \Delta^2 S_{n,k}
\bigr)^{-1}
Y_{k}^T \Delta s_n.
\label{eq:general_extrap}
\end{equation}
% where $\Delta S_{n,k}$ denotes the $p\times k$ matrix whose columns are
% \[
% \Delta s_n,\ \Delta s_{n+1},\ \ldots,\ \Delta s_{n+k-1},
% \]
% and $\Delta^2 S_{n,k}$ is the $p\times k$ matrix with columns
% \[
% \Delta^2 s_n,\ \Delta^2 s_{n+1},\ \ldots,\ \Delta^2 s_{n+k-1}.
% \]
% Here, the forward differences are defined by
% \[
% \Delta s_n = s_{n+1}-s_n,
% \qquad
% \Delta^2 s_n = \Delta s_{n+1}-\Delta s_n.
% \]
% Furthermore, $Y_{n,k}$ denotes a $p\times k$ matrix whose columns are given by
% \[
% y_n^{(1)},\ y_n^{(2)},\ \ldots,\ y_n^{(k)},
% \]
% and $Y_{n,k}^\ast$ stands for its conjugate transpose. 
We assume that all matrices involved in the inversion process are nonsingular which ensures the existence of the transformaed sequence. 
\\
In the particular case where $k=N$ and the vectors $y_i$ are chosen as the
canonical basis vectors $e_i$, $i=1,\ldots,N$, of $\mathbb{R}^N$, then the
transformation reduces to the {Henrici transformation} \cite[p.~116]{Henrici1964},
which was further studied in \cite{Brezinski1977}.
\medskip
The quantity $s_{k}^{(n)}$ can also be expressed as a ratio of two determinants.
Specifically,
\begin{equation}\label{s1}
s_{k}^{(n)}
=
\frac{
\begin{vmatrix}
s_n & s_{n+1} & \cdots & s_{n+k} \\
(y_1,\Delta s_n) & (y_1,\Delta s_{n+1}) & \cdots & (y_1,\Delta s_{n+k}) \\
\vdots & \vdots & & \vdots \\
(y_k,\Delta s_n) & (y_k,\Delta s_{n+1}) & \cdots & (y_k,\Delta s_{n+k})
\end{vmatrix}
}{
\begin{vmatrix}
1 & 1 & \cdots & 1 \\
(y_1,\Delta s_n) & (y_1,\Delta s_{n+1}) & \cdots & (y_1,\Delta s_{n+k}) \\
\vdots & \vdots & & \vdots \\
(y_k,\Delta s_n) & (y_k,\Delta s_{n+1}) & \cdots & (y_k,\Delta s_{n+k})
\end{vmatrix}
},
\end{equation}
where $(\cdot,\cdot)$ denotes the Euclidean inner product in $\mathbb{R}^N$ and where the determinant in the numerator is the vector obtained by expanding this determinant along the first row.. We also introduce the sequence $\beta_k^{(n)}$ obtained from the determinant expression of $s_{k}^{(n)}$ by replacing the first row in the numerator of \eqref{s1}) by: $\Delta s_n,\ldots, \Delta s_{n+k}$. The following algorithm that was called the $S\beta$-algorithm allows the recursive computation of the terms of the two sequences $s_k^{(n)}$ and $\beta_k^{(n)}$; see \cite{Jbilou1988,Jbilousadok1991} for more details.

\[
\left\{
\begin{array}{lll}
s_0^{(n)}=s_n,\; \beta_0^{(n)} = \Delta s_n, 
& n = 0,1,\ldots, \\[1.2ex]

\beta_k^{(n)} =
\displaystyle
\frac{\beta_{k-1}^{(n)} - a_k^{(n)}\,\beta_{k-1}^{(n+1)}}
     {1 - a_k^{(n)}}, 
& k \ge 1, \\[2ex]

s_k^{(n)} =
\displaystyle
\frac{s_{k-1}^{(n)} - a_k^{(n)}\, s_{k-1}^{(n+1)}}
     {1 - a_k^{(n)}}, \;\; 
a_k^{(n)} =
\displaystyle
\frac{\bigl( y_k , \beta_{k-1}^{(n)} \bigr)}
     {\bigl( y_k , \beta_{k-1}^{(n+1)} \bigr)}.
\end{array}
\right.
\]

\subsubsection{The vector H-algorithm} Let $(s_n)$ be a sequence of vectors and let $(g_i(n)$ some given sequences.  Define the new transformation ${\mathcal H}_k$ by $H_{k}^{(n)}={\mathcal H}_k(s_n)$ defined by the following ratio of two determinants
\begin{equation*}
H_{k}^{(n)}
=
\frac{
\begin{vmatrix}
s_n & s_{n+1} & \cdots & s_{n+k} \\
g_1(n) & g_1(n+1) & \cdots & g_1(n+k) \\
\vdots & \vdots & & \vdots \\
g_k(n) & g_k(n+1) & \cdots & g_k(n+k)
\end{vmatrix}
}{
\begin{vmatrix}
1 & 1 & \cdots & 1 \\
g_1(n) & g_1(n+1) & \cdots & g_1(n+k) \\
\vdots & \vdots & & \vdots \\
g_k(n) & g_k(n+1) & \cdots & g_k(n+k)
\end{vmatrix}
}.
\end{equation*}
As can be seen $s_k^{(n)}$ produced by the $S\beta$ algorithm is a particular case of $H_k^{(n)}$. The vectors $H_k^{n)}$ can
be recursively computed by the following vector  H-algorithm
\[
\left\{
\begin{array}{lll}
H_0^{(n)}=s_n,\; g_{0,i}^{(n)} = g_i(n), 
& n = 0,1,\ldots, \\[1.2ex]

H_k^{(n)} = H_{k-1}^{(n)}- 
\displaystyle \frac{g_{k-1,k}^{(n)}}{g_{k-1,k}^{(n+1)}-g_{k-1,k}^{(n)} }\, (H_{k-1}^{(n+1)}-H_{k-1}^{(n)}), 
& k \ge 1, \\[2ex]

g_{k,i}^{(n)} =g_{k-1,i}^{(n)}- 
\displaystyle \frac{g_{k-1,k}^{(n)}}{g_{k-1,k}^{(n+1)}-g_{k-1,k}^{(n)} }(g_{k-1,i}^{(n+1)}- g_{k-1,i}^{(n)}), 
&  \\[2ex]

\end{array}
\right.
\]
The vector H-algorithm was first introduced by Brezinski in \cite{Brezinski1983} and used to implement Henrici's transformation \cite{Henrici1964,Sadok1990}. Notice also that the H-algorithm is identical to the application of the scalar E-algorithm
simultaneously to each component of the vectors $s_n$  with the same auxiliary scalar sequences  $(g_i(n))$. A particular rule for jumping over breakdowns and near-breakdowns can also  be obtained for the H-algorithm, see \cite{Brezinski1991}.\\
Now, let  $\Delta S_{n+k-1}$ be the $N \times k$ matrix whose columns are $\Delta s_n,\ldots,\Delta s_{n+k-1}]$, ${\bf g}_n$ the $k$-vector whose elements are $g_1(n),\ldots,g_k(n)$ and $G_n$ the $k \times k$ matrix whose columns are $\Delta {\bf g}_n,\Delta {\bf g}_{n+1},\ldots,\Delta \Delta {\bf g}_{n+k-1}$. Then assuming that the matrix $G_n$ is non singular, it follows that
$$H_k^{(n)}= s_n-\Delta S_{n+k-1}\, (G_n)^{-1}\, {\bf g}_n,$$
which can also be expressed as the following Schur complement
\begin{equation*}
  H_k^{(n)} = \left(  {\mathcal H}_k^{(n)} / G_n \right), 
\end{equation*}
where ${\mathcal H}_k^{(n)}$ is the matrix given by
\begin{equation*}
{\mathcal H}_k^{(n)}=
\begin{pmatrix}
s_n & \Delta S_{n+k-1}\\
{\bf g}_n & G_n
\end{pmatrix}.
\end{equation*}

\subsubsection{Application of VEP to linear systems of equations}

Consider the linear system
\begin{equation}\label{eqe1}
C x = f,
\end{equation}
with $C \in \mathbb{R}^{N \times N}$ and $f \in \mathbb{R}^N$. Let $M$ be a nonsingular preconditioner. The preconditioned system is
\begin{equation}
M^{-1} C x = M^{-1} f.
\end{equation}
Starting from $s_0$, generate
\begin{equation}
s_{j+1} = B s_j + b, \quad j=0,1,\ldots,
\end{equation}
with $B = I - A$, $A = M^{-1} C$, and $b = M^{-1} f$. Then the generalized residual becomes the true residual:
\begin{equation}
\tilde r(t_k^{(n)}) = r(t_k^{(n)}) = b - A t_k^{(n)}.
\end{equation}
Moreover, $\Delta^2 s_n = -A \Delta s_n$ implies $\Delta^2 S_{k,n} = -A \Delta S_{k,n}$. Next, and for simplification, we set $n=0$ and denote $t_k=t_k^{(0)}$. 
Let $d$ be the degree of the minimal polynomial of $B$ with respect to $s_0 - x^*$. Then for $k \le d$, $\Delta S_k$ and $\Delta^2 S_k$ are of full rank, and $t_d$ exists and equals the solution of equation \eqref{eqe1}. 
The residuals generated by RRE, MPE, and MMPE satisfy
$$
\left\{
\begin{array}{ll}
r_k \in \widetilde W_k = A \widetilde V_k,\\
r_k \perp \widetilde L_k,
\end{array}
\right.
$$
where $\widetilde V_k = \mathrm{span}\{\Delta s_0, \dots, \Delta s_{k-1}\}$ and
\[
\widetilde L_k =
\begin{cases}
\widetilde W_k, & \text{RRE}, \\
\widetilde V_k, & \text{MPE}, \\
\widetilde Y_k = \mathrm{span}\{y_1, \dots, y_k\}, & \text{MMPE}.
\end{cases}
\]
Let $P_k$ be the orthogonal projector onto $\widetilde W_k$. Then the RRE residual is given by 
\begin{equation}\label{eqrk}
r_k^{\text{RRE}} = r_0 - P_k r_0.
\end{equation}
Let $Q_k$ and $R_k$ denote oblique projectors for MPE and MMPE, yielding
\begin{align*}
r_k^{\text{MPE}} &= r_0 - Q_k r_0,\; \text{and}\\
r_k^{\text{MMPE}} &= r_0 - R_k r_0.
\end{align*}
Define the acute angle $\theta_k$ between $r_0$ and the orthogonal projection $P_kr_0$
\begin{equation}
\cos \theta_k = \max_{z \in \widetilde W_k \setminus \{0\}} \frac{|(r_0, z)|}{\|r_0\| \, \|z\|}. 
\end{equation}
We also define the two acuate angles  $\phi_k$ and $\psi_k$ between the vectors $r_0, Q_kr_0$ and $r_0,\R_kr_0$, respectively.\\
In the following, we give some relations satisfied by the residual norms of the three
extrapolation methods.

\medskip
\begin{theorem}\cite{Jbilousadok2000}
Let $\theta_k$ be the acuate angle between  $r_0$ and the orthogonal projection $P_kr_0$, let $\phi_k$ be the acute angle between $r_0$ and $Q_k \, r_0$ and let $\psi_k$
denote the acute angle between $r_0$ and $R_k \, r_0$. Then we have the following relations\\
\begin{enumerate}
\item $\parallel r_k^{rre} \parallel^2=(1-\cos^2 \, \theta_k)\,
\parallel r_0 \parallel^2$.
\item $\parallel r_k^{mpe} \parallel^2=(tan^2 \, \phi_k)\,
\parallel r_0 \parallel^2$.
\item $\parallel r_k^{mpe} \parallel\, \le \, (\cos \, \phi_k)\,
\parallel r_k^{rre} \parallel$.
\end{enumerate}

\medskip
\nd Moreover if for MMPE $y_j=r_0$ for some $j=1,\ldots,k$, then we also have

\begin{itemize}
\item $\parallel r_k^{mmpe} \parallel^2=(tan^2 \, \psi_k)\,
\parallel r_0 \parallel^2$.
\item  $\parallel r_k^{mmpe} \parallel\, \le \, (\cos \, \psi_k)\,
\parallel r_k^{rre} \parallel$.
\end{itemize}
\end{theorem}

\medskip

\nd The general vector extrapolation algorithm for solving linear systems of equations is summarized in Algorithm \ref{VEPL}.

\begin{algorithm}[h]
\caption{Vector polynomial extrapolation (VPE) for linear systems}
\label{VEPL}
Initialize $k = 0$, choose $x_0$, and set integers $p$ and $m$.
\begin{enumerate}
    \item {Basic iteration:} Set $t_0 = x_0$, $z_0 = t_0$, then for $j = 0, \dots, p-1$ compute
    \[
    z_{j+1} = B z_j + b.
    \]
    \item {Extrapolation:} Set $s_0 = z_p$, then for $j = 0, \dots, m$ compute
    \[
    s_{j+1} = B s_j + b,
    \]
    and compute the approximation $t_m$ using RRE, MPE, or MMPE.
    \item Update $s_0 = t_m$, increment $k = k + 1$, and return to step 1 until convergence.
\end{enumerate}
\end{algorithm}
\medskip

\nd Stable implementations of these methods use either the QR decomposition \cite{Sidi1991} for RRE and MPE  or LU decomposition with pivoting strategies  \cite{Jbilousadok1999} for MMPE,  to reduce storage and improve numerical stability.\\ 
A key property of vector polynomial extrapolation methods is the fact that they  are mathematically equivalent to some well known Krylov subspace methods such as GMRES \cite{SaadS1986}  and Full Orthogonalisation  Method (FOM) \cite{Saad2003} when they are applied to vector sequences generated linearly ; see \cite{Beneu1983,Jbilou1988,Jbilousadok2000} for more details. This is stated in the following theorem.

\medskip
\begin{theorem}
    When applied to linearly generated sequences, the vector extrapolation methods RRE and MPE are mathematically equivalent to GMRES and Arnoldi methods, respectively while  MMPE is equivalent to an oblique projection Krylov subspace method. 
\end{theorem}
\medskip

\nd From the results of this theorem, we can state   that all polynomial vector extrapolation (PVE) methods have  an interpretation in terms of polynomial approximation theory. The
coefficients $\{\alpha_j\}$ determine a polynomial $P(\lambda)$ whose values are small at the dominant eigenvalues
$\lambda_j$ of the iteration matrix, that is  $
P(\lambda_j) \approx 0$.
As a consequence, the resulting extrapolated iterates effectively suppress the principal components of the error, thereby enhancing the convergence rate. \\
Notice also that since the VPE methods are mathematically equivalent to Krylov subspace methods when applied for solving linear systems, they converge in a maximum $N$ iterations where $N$ is the dimension of the space. However and as for Krylov subspace methods, extrapolation methods are usually used in a restarted mode and this allow saving storage for large problems.

\subsubsection{Application to linear-discrete ill-posed problems}\label{subsubsec4}
In this subsection, we focus exclusively on applying the Reduced Rank
Extrapolation (RRE) to the sequence of vectors generated by the Truncated
Singular Value Decomposition (TSVD). The TSVD iterates are viewed as an ordered
sequence of approximations, where the truncation index plays the role of an
iteration parameter. In this framework, the inclusion of additional singular
components progressively enriches the approximation, but may also introduce
instability due to noise amplification in the case of ill-posed problems.
\\
The purpose of applying RRE is to accelerate the convergence of the TSVD
sequence while preserving its regularizing properties. By forming appropriate
linear combinations of consecutive TSVD iterates, RRE exploits the underlying
low-rank structure of the problem and enhances the contribution of the most
relevant singular modes. Consequently, RRE acts as an effective post-processing
technique for TSVD-based regularization.
\\
Furthermore, the RRE approach allows for an implicit and adaptive control of
the effective rank of the approximation. The resulting extrapolated solutions
often achieve a favorable compromise between accuracy and stability, making
the combined TSVD--RRE strategy particularly attractive for large-scale
ill-posed problems. For more details and explanations, we refer to the paper \cite{Jbilousadokreichel2009}.
\\

We consider the problem of approximating a solution to the linear system
\begin{equation}\label{linsys}
Ax = b,
\end{equation}
where $A \in \mathbb{R}^{m \times n}$ is a matrix with ill-determined rank; that is, $A$ has many singular values of widely varying magnitudes, some of which are close to zero. In particular, $A$ is severely ill-conditioned and may be singular. The method is described here for the case $m \geq n$, although it is also applicable when $m < n$.
\\
Let $\hat{b} \in \mathbb{R}^m$ denote the unknown exact right-hand side vector associated with $b$, so that
\begin{equation}\label{error}
b = \hat{b} + e,
\end{equation}
where $e$ represents the measurement or computational error. Our goal is to compute an accurate approximation of the minimal-norm solution $\hat{x}$ of the exact consistent system
\[
A x = \hat{b}.
\]
Formally, $\hat{x} = A^\dag \hat{b}$, where $A^\dag$ denotes the Moore-Penrose pseudoinverse of $A$. Due to the presence of the error $e$ and the ill-conditioning of $A$, the direct application of $A^\dag$ to $b$,
\[
A^\dag b = A^\dag (\hat{b} + e) = \hat{x} + A^\dag e,
\]
generally does not yield a reliable approximation of $\hat{x}$. 
A standard approach for approximating $\hat{x}$ is to employ the Truncated Singular Value Decomposition (TSVD), which replaces $A^\dag$ with a low-rank approximation; see, e.g., Golub and Van Loan \cite{GolubVL1996} or Hansen \cite{Hansen1994}. Let the singular value decomposition of $A$ be given by
\begin{equation}\label{svdA}
A = \sum_{j=1}^n \sigma_j u_j v_j^T,
\end{equation}
where the singular values satisfy
\begin{equation}\label{singval}
\sigma_1 \geq \sigma_2 \geq \cdots \geq \sigma_\ell > \sigma_{\ell+1} = \cdots = \sigma_n = 0.
\end{equation}
Here, the vectors $u_j \in \mathbb{R}^m$ and $v_j \in \mathbb{R}^n$ form the columns of the orthonormal matrices $U = [u_1, \dots, u_n] \in \mathbb{R}^{m \times n}$ and $V = [v_1, \dots, v_n] \in \mathbb{R}^{n \times n}$. For the problems considered, many of the smallest nonzero singular values are very close to zero.
\\
For any $1 \leq k \leq \ell$, the rank-$k$ approximation of $A$ is defined as
\begin{equation}\label{Ak}
A_k = \sum_{j=1}^k \sigma_j u_j v_j^T,
\end{equation}
with corresponding pseudoinverse
\begin{equation}\label{Akinv}
A_k^\dag = \sum_{j=1}^k \sigma_j^{-1} v_j u_j^T.
\end{equation}
The minimal-norm solution of the least-squares problem
\[
\min_{x \in \mathbb{R}^n} \|A_k x - b\|
\]
is then given by
\begin{equation}\label{xk}
x_k := A_k^\dag b = \sum_{j=1}^k \frac{u_j^T b}{\sigma_j} v_j.
\end{equation}
A key question is the selection of an index $k$ such that $x_k$ provides an accurate approximation of $\hat{x}$. To formalize this, let $k_{\rm opt} \geq 1$ denote the smallest integer satisfying
\begin{equation}\label{kopt}
\|x_{k_{\rm opt}} - \hat{x}\| = \min_{k \geq 1} \|x_k - \hat{x}\|.
\end{equation}

\nd In the sequel, we set
\[
y_i=\Delta^2 s_{i-1};\;\;\mbox{for }\,
1\leq i\leq m-2,
\]
and consider the sequence $\{s_k\}_{k\geq 0}$ generated by TSVD. Thus,
\begin{equation}\label{sk}
s_0:=0,\qquad s_k:=x_k=A_k^\dagger b=\sum_{j=1}^k
\frac{u_j^Tb}{\sigma_j}v_j=\sum_{j=1}^k {\delta_j}v_j,
\end{equation}
where
\begin{equation}\label{deltajdef}
\delta_j:=\frac{u_j^Tb}{\sigma_j},\qquad 1\leq j\leq k.
\end{equation}
Thus, we have
\begin{equation}\label{formule0}
 \Delta
s_{k-1}=s_k-s_{k-1}={\delta_k}v_k,
\end{equation}
We may assume that $\delta_k \ne 0$, because if this is not the case, then we
delete the corresponding member from the sequence (\ref{sk}) and
compute the next one by keeping the same index notation. 
The
matrix $\Delta S_{k-1}=[\Delta s_{0},\ldots,\Delta s_{k-1}]$ can
be factored according to
\begin{equation}\label{formule2}
\Delta
S_{k-1}=\left[{\delta_1}v_1,\ldots,\delta_{k}v_{k}\right]=V_k\;
\mbox{diag}[{\delta_1},\ldots,{\delta_{k}}],
\end{equation}
where $V_k=[v_1,\ldots,v_{k}]$. Moreover, since
$\Delta^2 s_{k-1}={\delta_{k+1}}v_{k+1}-{\delta_k}v_k,$ we deduce that
\begin{equation}\label{formule1}
 \Delta^2 S_{k-1} =V_{k+1}\; \left[\begin{array}{cccc}
  -{\delta_1} &  & &  \\
{\delta_2} & -{\delta_2} & &  \\
    &\ddots &\ddots &  \\
   &   & {\delta_k}&  -{\delta_k}\\
   & & & {\delta_{k+1}} \
 \end{array}\right].
 \end{equation}
 Then using (\ref{formule2}), we get  $$\Delta^2 S_{k-1}^T \Delta
S_{k-1}= \left[\begin{array}{cccc}
  -{\delta_1} &  & &  \\
{\delta_2} & -{\delta_2} & &  \\
    &\ddots &\ddots &  \\
   &   & {\delta_k}&  -{\delta_k}\\
   & & & {\delta_{k+1}} \
 \end{array}\right]^T V_{k+1}^T \;
 [\delta_1\, v_1,\ldots,\delta_k\,v_k].$$
 On the other hand, since
\begin{equation}\label{formule3}
 V_{k+1}^T\,
 \left[{\delta_1}v_1,\ldots,{\delta_{k}}v_{k}\right]=
 \left[
 \begin{array}{cccc}
 \delta_1&0&\ldots &0\\
0&\delta_2&\ldots &0\\
\vdots&\vdots&\ddots & \vdots\\
0&\ldots & 0 & \delta_k\\
0&0&\ldots &0
\end{array}
\right],
\end{equation}
it follows that
 $$
\Delta^2 S_{k-1}^T \Delta S_{k-1}=\left[\begin{array}{cccc}
  -\delta_1^2 &\delta_2^2 & &  \\
  &- \delta_2^2 &\delta_3^2 &  \\
 &&&\\
 &&&\\
   &   & -\delta_{k-1}^2& \delta_k^2\\
   & & &- \delta_{k}^2 \
 \end{array}\right].$$
The extrapolated vector $t_k=t_k^{(0)}$ is given by
\begin{equation}\label{tk1}
t_k=\displaystyle \sum_{j=0}^k \gamma_j^{(k)} s_{j},
\end{equation}
where the vector $\boldsymbol{\gamma}$ whose components are $\gamma_0^{(k)},\ldots,\gamma_k^{(k)}$ solves the following linear system of equations
\begin{equation}\label{lr1}
\displaystyle \sum_{j=0}^k \gamma_j^{(k)}=1,\;\; \;\; \displaystyle \sum_{j=0}^k \gamma_j^{(k)} \eta_{ij}=0,
\end{equation}
where $\eta_{ij}=(\Delta^2 s_i,\Delta s_j)$.\\
Let $\beta_i = \gamma_i^{(k)} / \gamma_k^{(k)}$, then each coefficient $\gamma_i^{(k)}$ can be recovered from the normalized variables $\beta_i$ through
\[
\gamma_i^{(k)} = \frac{\beta_i^{(k)}}{\displaystyle \sum_{j=0}^k \beta_j^{(k)}}.
\]
It follows that the vector
\[
\boldsymbol{\beta}_k = (\beta_0^{(k)},\ldots,\beta_{k-1}^{(k)})^{T}
\]
satisfies the following linear system of equations:
\begin{equation}\label{lin2}
\left(\Delta^2 S_{k-1}^{T}\,\Delta S_{k-1} \right)\,\boldsymbol{\beta}_k
=
-\,\Delta^2 S_{k-1}^{T}\,\Delta s_k.
\end{equation}
Due to the structure of the finite difference matrices, the right-hand side simplifies to
\[
\Delta^2 S_{k-1}^{T}\,\Delta s_k
=
[0,\ldots,0,\delta_{k+1}^2]^{T},
\]
which immediately yields the explicit solution
\[
\beta_i^{(k)} = \frac{\delta_{k+1}^2}{\delta_{i+1}^2},
\qquad 0 \le i < k.
\]
Summing these coefficients, we obtain
\[
\sum_{i=0}^k \beta_i^{(k)}
=
1 + \sum_{i=0}^{k-1} \frac{\delta_{k+1}^2}{\delta_{i+1}^2}
=
\delta_{k+1}^2 \sum_{i=0}^k \frac{1}{\delta_{i+1}^2}.
\]
Consequently, the coefficients $\gamma_j^{(k)}$ admit the closed-form expression
\begin{equation}\label{gama}
\gamma_j^{(k)}
=
\frac{\displaystyle \frac{1}{\delta_{j+1}^2}}
     {\displaystyle \sum_{i=0}^k \frac{1}{\delta_{i+1}^2}},
\qquad 0 \le j < k.
\end{equation}
Assume now that the coefficients $\gamma_0^{(k)},\gamma_1^{(k)},\ldots,\gamma_k^{(k)}$ have been computed. We introduce the auxiliary variables
\begin{equation}
\alpha_0^{(k)} = 1 - \gamma_0^{(k)},
\qquad
\alpha_j^{(k)} = \alpha_{j-1}^{(k)} - \gamma_j^{(k)},\quad 1 \le j < k,
\qquad
\alpha_{k-1}^{(k)} = \gamma_k^{(k)},
\label{eq:alpha-def}
\end{equation}
which allow us to rewrite the extrapolated vector $t_k$ in a telescopic form:
\begin{equation}\label{ttk}
t_k
=
s_0 + \sum_{j=0}^{k-1} \alpha_j^{(k)} \,\Delta s_j
=
s_0 + \Delta S_{k-1}\,\boldsymbol{\alpha}_k,
\end{equation}
where $ 
\boldsymbol{\alpha}_k = [\alpha_0^{(k)},\ldots,\alpha_{k-1}^{(k)}]^{T}$. 
Using \eqref{gama}, the coefficients $\alpha_i$ can be expressed explicitly as
\begin{equation}\label{alpha}
\alpha_i^{(k)}
=
\sum_{j=i+1}^{k} \gamma_j^{(k)}
=
\frac{\displaystyle \sum_{j=i+1}^{k} \frac{1}{\delta_{j+1}^2}}
     {\displaystyle \sum_{l=0}^{k} \frac{1}{\delta_{l+1}^2}},
\qquad 0 \le i < k.
\end{equation}
Finally, the extrapolated vector $t_k$ admits the spectral representation
\begin{equation}\label{tk}
t_k
=
\sum_{j=1}^{k}
\alpha_{j-1}^{(k)}
\frac{u_j^{T} b}{\sigma_j}\, v_j.
\end{equation}
This representation reveals that the application of the RRE method acts as a filtering mechanism for the truncated singular value decomposition (TSVD), with filter factors
\[
f_j^{(k)} = \alpha_{j-1}^{(k)},
\qquad 1 \le j \le k.
\]
We need to evaluate $\Vert t_{k+1}-t_k\Vert$ for our stopping criterion.
From (\ref{tk}) we obtain that
\begin{equation*}
t_{k+1}-t_k=\sum_{j=1}^k \frac{( \alpha_{j-1}^{(k+1)}- \alpha_{j-1}^{(k)})}
{\sqrt{\delta_j}} v_j+\dfrac{\alpha_{k}^{(k+1)}}{\sqrt{\delta_{k+1}}} v_{k+1}.
\end{equation*}
Since the vectors $v_j$, $1\leq j\leq k+1$, are orthonormal, it follows that
\begin{equation*}
\Vert t_{k+1}-t_k\Vert=\sqrt{\sum_{j=1}^k
\frac{\vert \alpha_{j-1}^{(k+1)}- \alpha_{j-1}^{(k)}\vert^2}
{\delta_j}+\dfrac{\vert\alpha_{k}^{(k+1)}\vert^2}{\delta_{k+1}}}.
\end{equation*}
The generalized residual can be written as
\begin{equation}
\tilde r(t_k)=\Delta S_k\, \boldsymbol{\gamma}_k=V_{k+1}\;
\mbox{diag}\left[\frac{u_1^Tb}{\sigma_1},\ldots,
\frac{u_{k+1}^Tb}{\sigma_{k+1}}\right] \boldsymbol{\gamma}_k,
\end{equation}
where
$\boldsymbol{\gamma}_k=[\gamma_0^{(k)},\gamma_1^{(k)},\ldots,\gamma_k^{(k)}]^T$
is given by (\ref{gama}). It follows from (\ref{deltajdef}) that
$$\Vert \tilde r(t_k)\Vert^2=\displaystyle
\sum_{i=0}^k {\delta_{i+1}^2\,\left(\gamma_i^{(k)} \right)^2},  $$
and by using (\ref{gama}), we obtain the simple expression
$$\Vert \tilde r(t_k)\Vert=\displaystyle \frac{1}{\sqrt{\displaystyle
\sum_{j=0}^{k} \frac{1}{\delta_{j+1}^2}}}.$$
Moreover, from the relation (\ref{gama}) we also have an expression of the norm of the generalized residual:
$$\Vert \tilde r(t_k)\Vert = \sqrt{ \delta_{k}^2 \gamma_{k-1}^{(k)}}. $$
For well-posed problems, the norm $\lVert \tilde r(t_k) \rVert$ decreases
monotonically and tends to zero as $k$ increases. In contrast, for ill-posed
problems affected by perturbations in the right-hand side, the quantity
$\lVert \tilde r(t_k) \rVert$ initially decreases as $k$ increases, provided
that $k$ remains sufficiently small. Beyond a certain threshold, however, the
residual norm begins to increase with increasing $k$.
The value of the index $k$ at which
$\lVert \tilde r(t_k) \rVert$ ceases to decrease often provides an approximation
$x_k$ that is close to the desired solution $\hat{x}$.\\
The RRE-TSVD algorithm is summarized as follows:\\
\begin{algorithm}[h]
 \caption{The RRE-TSVD
  algorithm for ill-posed problems}
\begin{itemize}
\item {Compute the SVD of the matrix $A$}:\;
$[U,\Sigma,V]=svd(A)$.\\ { Set} $s_0=0$, $s_1=\frac{u_1^Tb}{\sigma_1}v_1$,
{ and} $t_1=s_1$, { with} \\
$u_i=U(:,i)$ { and} $v_i=V(:,i)$ { for} $i=1,\ldots,n$.
 \item {For} $k=2,\ldots, n$
\begin{enumerate}
 \item {Compute } $s_k$ {from } (\ref{sk}).
 \item {Compute the} $\gamma_i^{(k)}$ { and} $\alpha_i^{(k)}$ {for}
 $i=0,\ldots,k-1$, { using} (\ref{gama}) {\tt and} (\ref{alpha}).
 \item { Form the approximation} $t_k$ { by} (\ref{tk}).
 \item {If} $\Vert t_{k}-t_{k-1} \Vert / \Vert t_{k-1} \Vert <
 tol$, { stop}.
 \end{enumerate}
 \item End
\end{itemize}
 \end{algorithm}

 \nd We illustrate the convergence behavior of the
sequences $\{t_k\}_{k\geq 1}$ and $\{\Vert \tilde r(t_k)\Vert\}_{k\geq 1}$ by
considering the
integral equation
\begin{equation}\label{bart}
\int_{0}^{\pi/2} \kappa(s,t)x(t)dt=g(s),\qquad 0 \le s \le \pi,
\end{equation}
where$$ \kappa(s,t)=\exp(s \, \cos(t))\;\; {\rm and} \;\; g(s) = 2
\sin(s)/s.$$ The solution is given by   $x(t) = \sin(t)$. This
integral equation is discussed in \cite{Baart1982}. We used the MATLAB
code {\sf baart} from \cite{Hansen1998} to discretize (\ref{bart}) by a
Galerkin method with $600$ orthonormal box functions as test and trial
functions. This yields the nonsymmetric matrix $A \in {\R}^{600\times 600}$
with condition number $\kappa(A):=\Vert A
\Vert \, \Vert A^{-1} \Vert =4\cdot 10^{18}$. Thus, $A$ is numerically
singular. The vectors $\hat{b}$ and $b$ are determined as described above with
noise-level $\nu=1\cdot 10^{-2}$.

 \begin{figure}[h] 
 \centering
\includegraphics[width=0.4\textwidth]{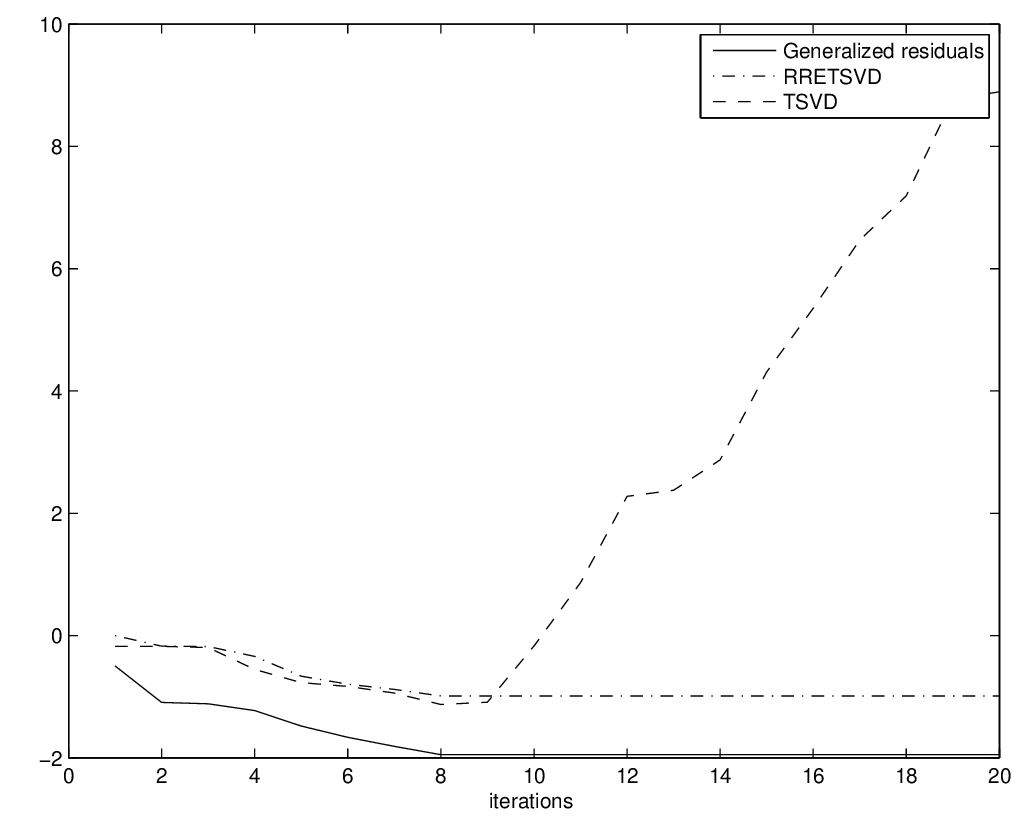}
    \caption{ The relative errors}\label{fig1}
\end{figure}

% % $\Vert s_k-\hat{x}\Vert/\Vert \hat{x} \Vert$ for TSVD (dashed graph),
% % the relative errors $\Vert t_k-\hat{x}\Vert/\Vert \hat{x} \Vert$ for RRE-TSVD
% % (dash-dotted graph), and the norm of the generalized residuals
% % $\Vert \tilde r(t_k) \Vert$ (solid graph).}
% % \label{fig1}

\nd Figure \ref{fig1} shows the convergence of TSVD and RRE-TSVD.
Observe from  this figure that when $k$ increases and is
sufficiently small, the relative error norms
\[
\|s_k-\hat{x}\|/\|\hat{x}\|~~\mbox{and}~~\|t_k-\hat{x}\|/\|\hat{x}\|
\]
The sequence $(s_k)$ is generated by the TSVD method, whereas the sequence
$(t_k)$ is obtained by applying RRE to the TSVD iterates. For small values of
$k$, the approximation errors in both sequences decrease steadily. When the
truncation index exceeds the optimal value $k_{\rm opt}=9$, the error
associated with $s_k$ begins to increase due to noise amplification, while
the error norm of $t_k$ remains essentially constant and attains its minimum
at $k=k_{\rm opt}$. A comparable stagnation behavior is observed in the
generalized residual norm.
\\
This observation indicates that the point at which the generalized residual
norm ceases to decrease provides a reliable and computationally inexpensive
criterion for identifying the optimal truncation index. As a result,
monitoring this stagnation offers a practical strategy for choosing $k$ in
both the TSVD and RRE--TSVD approaches (see equation~\eqref{kopt} for the
definition of $k_{\rm opt}$).

\subsubsection{Application of VEP to nonlinear problems}

Consider the nonlinear system
\begin{equation}
G(x) = x, \quad G:\mathbb{R}^N \to \mathbb{R}^N,
\end{equation}
with residual $r(x) = G(x) - x$. Starting from $s_0$, generate the sequence
\begin{equation}
s_{j+1} = G(s_j), \quad j = 0,1,\dots,
\end{equation}
 so that $r(s_j) =  \Delta s_j$.
%Preconditioning may be applied to improve convergence:
% \begin{equation}
% \tilde G(x) = x.
% \end{equation}
The extrapolation algorithm is summarized as follows:
\begin{algorithm}[H]
\caption{Vector Extrapolation Method for Nonlinear Systems}\label{VEPN}
\begin{enumerate}
    \item Set $k=0$, choose initial guess $x_0$, integers $p$ and $m$, and tolerance $\epsilon$.
    \item \textbf{Basic iteration:} 
    \begin{enumerate}
        \item Set $t_0 = x_0$, $w_0 = t_0$.
        \item For $j = 0, \dots, p-1$, compute $w_{j+1} =  G(w_j)$.
    \end{enumerate}
    \item \textbf{Extrapolation phase:} 
    \begin{enumerate}
        \item Set $s_0 = w_p$.
        \item If $\|s_1 - s_0\| < \epsilon$, stop.
        \item Otherwise, for $j = 0, \dots, m$, compute $s_{j+1} =  G(s_j)$.
        \item Compute $t_m$ using RRE, MPE, or MMPE.
    \end{enumerate}
    \item Set $s_0 = t_m$, $k = k+1$, $x_0 = s_0$, and go back to step 2.
\end{enumerate}
\end{algorithm}

\medskip
\nd  
Let us denote 
\[
G'(x^*) = J,
\] 
where $x^*$ is a solution and assume  the function \(G\) satisfies the following conditions.
\\
Notice that these vector  methods do not require the Jacobian of $\tilde G$ and, in complete form, achieve quadratic convergence under some conditions as stated in the next theorem.\\
The matrix \(J - I\) is nonsingular. We set
\begin{equation}
M = (J - I)^{-1}. \label{eq:M_def}
\end{equation}
Furthermore, the Fréchet derivative \(G'\) of \(G\) satisfies a Lipschitz condition, i.e.,
\begin{equation}
\| G'(x) - G'(y) \| \le L \, \| x - y \|, \quad \text{for all } x, y \in D,
\end{equation}
where \(L\) is a positive constant and \(D\) denotes the domain of interest.
\medskip

\begin{theorem}\cite{Jbilousadok1991}
Under some more conditions on \(G\) defined in \cite{Jbilousadok1991}, there exists a neighborhood \(U\) of \(x^*\) such that, for all
\(x_0 \in U\),
\[
\|x_{k+1} - x^*\| = o\!\left(\|x_k - x^*\|^2\right),
\]
for both  RRE and MPE methods.
\end{theorem}
\medskip

 In practice, these extrapolation techniques are applied in a black-box manner to the
sequence of iterates $(s_k)$ generated by a nonlinear process
\[
s_{i+1} = G(s_i), \; i \ge 0,
\]
without requiring explicit access to the Jacobian $G'(x^\ast)$. When the nonlinear
iteration converges locally and the mapping $G$ is sufficiently smooth, the associated
error sequence often admits an asymptotic representation as a finite sum of linear modes.
Under these circumstances, extrapolation methods are capable of canceling the dominant
error components and thereby accelerating convergence.
\\
Among the three approaches, RRE is generally regarded as the most robust in nonlinear
settings. Its formulation as a least-squares problem provides increased stability with
respect to perturbations, roundoff errors, and deviations from the linear model. For this
reason, RRE is frequently preferred in large-scale nonlinear applications, including the
numerical solution of nonlinear partial differential equations, integral equations, and
inverse problems.
\\
The MPE method, while optimal for linear problems in a polynomial sense, may display less
predictable behavior for nonlinear iterations. In particular, the linear error model
underlying MPE can be inaccurate away from the fixed point, which may result in stagnation
or even divergence. Nevertheless, when the iterates are sufficiently close to the
solution, MPE can still provide substantial acceleration and, in some cases, outperform
RRE.\\
The MMPE method can be viewed as an intermediate approach between RRE and MPE. By
introducing additional flexibility through auxiliary vectors or weights, MMPE may improve
stability in certain nonlinear problems. However, its effectiveness depends more strongly
on parameter selection and problem structure, which limits its widespread use in general
nonlinear applications. 
The performance of RRE, MPE, and MMPE in nonlinear problems relies on the validity
of a local linearization of the iteration near the solution. When this approximation is
adequate, extrapolation methods can significantly enhance convergence. Otherwise,
safeguard techniques such as restarting, truncation, or adaptive window strategies are
often required to maintain stability.

\subsection{The $\epsilon$-algorithms} In this part, we consider two well known vector $\varepsilon$-algorithms: the vector $\varepsilon$ of Wynn and the  topological $\varepsilon$-algorithm of Brezinski, and give also  some variants of these algorithms.

\subsubsection{The vector $\epsilon$-Algorithm}

The scalar $\epsilon$-algorithm, introduced by Wynn \cite{Wynn1962}, is a widely used method for accelerating the convergence of scalar sequences. A natural question is how to generalize this algorithm to vector sequences $(s_n) \subset \mathbb{R}^N$ or $\mathbb{C}^N$. A main difficulty in this generalization arises from the fact that vectors do not have a natural multiplicative inverse as scalars do. 
\\
One approach, proposed by Wynn \cite{Wynn1962}, is to define the inverse of a vector $z \in \mathbb{R}^N$ as
\begin{equation}
z^{-1} = \frac{z}{\|z\|^2},
\end{equation}
where $\|\cdot\|$ denotes the Euclidean norm in $\mathbb{R}^N$. This definition preserves certain algebraic properties needed to extend the scalar algorithm to vectors.

\medskip
\nd With this vector inverse, the \emph{vector $\epsilon$-algorithm} (VEA) can be defined recursively as
\begin{equation}
\left\{
\begin{array}{l}
\epsilon_{-1}^{(n)} = 0, \quad \epsilon_0^{(n)} = s_n, \quad n = 0,1, \ldots,\\[1em]
\epsilon_{k+1}^{(n)} = \epsilon_{k-1}^{(n+1)} + \big[\epsilon_k^{(n+1)} - \epsilon_k^{(n)}\big]^{-1}, \quad k,n = 0,1,\ldots
\end{array}
\right.
\end{equation}

\medskip
\noindent
The convergence properties of this algorithm have been studied extensively. It was proved in \cite{Graves1992,Mcleod1971} that if there exists a finite set of scalars $a_0, \dots, a_k$ such that if 
\begin{equation}
\sum_{i=0}^k a_i (s_{n+i} - s) = 0, \quad a_k \neq 0, \quad a_0 + \dots + a_k \neq 0,
\end{equation}
for all $n \ge N_0$, then the vector $\epsilon$-algorithm satisfies
\begin{equation}
\epsilon_{2k}^{(n)} = s, \quad \forall n \ge N_0.
\end{equation}

\medskip
\nd The vector $\epsilon$-algorithm can also be applied to sequences arising from the iterative solution of linear systems. Specifically, if $(s_n)$ is generated by an iterative scheme of the form
\begin{equation}
s_{n+1} = B s_n + b,
\end{equation}
then both the scalar and vector $\epsilon$-algorithms yield the exact solution $s^*$ of the linear system after a finite number of steps:
\begin{equation}
\epsilon_{2N}^{(n)} = s^*, \quad \forall n \ge 0,
\end{equation}
The intermediate quantities $\epsilon_{2k}^{(n)}$, with $k < N$, provide increasingly accurate approximations of $x^*$.

\medskip
\nd The vector $\epsilon$-algorithm has also been successfully applied to nonlinear sequences, allowing acceleration of fixed-point iterations for nonlinear problems. For instance, it can be applied to sequences generated by a nonlinear operator $G$ via
\begin{equation}
s_{n+1} = G(s_n),
\end{equation}
to accelerate convergence toward a fixed point $s^* = G(s^*)$.

\medskip
\noindent
 Although powerful, the vector $\epsilon$-algorithm is more computationally demanding than vector polynomial extrapolation methods such as RRE, MPE, and MMPE. Specifically:
\begin{itemize}
\item Computing $\epsilon_{2k}^{(n)}$ requires $2k+1$ sequence vectors, $s_n, \dots, s_{n+2k}$, whereas the polynomial methods only require $k+2$ vectors, $s_n, \dots, s_{n+k+1}$.
\item Each step involves vector inversions and additional storage for the intermediate $\epsilon$-quantities, which can be significant for large $n$ or high $k$.
\item Despite the higher storage and computational requirements, the algorithm provides highly accurate approximations and is particularly useful when sequences exhibit slow convergence or oscillatory behavior.
\end{itemize}

\medskip
\nd In summary, the vector $\epsilon$-algorithm is a powerful convergence acceleration tool for both linear and nonlinear sequences, complementing the vector polynomial methods.  Its main advantage lies in its strong convergence properties, while its main limitation is the higher computational effort required. More properties and developments on the vector $\epsilon$-algorithm can be found in \cite{Brezinski1991,GravesRS2000,Salam1996}. 

A later and more sophisticated development was based on the use of
noncommutative algebraic structures. In this setting, the classical concept of
determinants is replaced by \emph{designants}, which provide an appropriate
generalization in noncommutative fields. This approach was further developed
within the framework of real Clifford algebras \cite{Salam1996}. Both left and
right versions of designants were defined and analyzed in early foundational
work \cite{Heyting1927}.
\\
To clarify the construction, consider the right designant
$\Delta^{(n)}_r$ associated with the matrix
\[
\Delta^{(n)}_r =
\begin{vmatrix}
a_{11} & \cdots & a_{1n} \\
\vdots &        & \vdots \\
a_{n1} & \cdots & a_{nn}
\end{vmatrix}_r,
\]
where the coefficients $a_{ij}$ take their values in a noncommutative field.
Rather than being defined by permutation expansions, this quantity is introduced
through a recursive procedure. In the simplest nontrivial case, namely for a
$2\times2$ matrix, one obtains
\[
\Delta^{(2)}_r =
\begin{vmatrix}
a_{11} & a_{12} \\
a_{21} & a_{22}
\end{vmatrix}_r
= a_{22} - a_{12} a_{11}^{-1} a_{21}.
\]
In the general case, let $A^{r}_{pq}$ denote the right designant of order $p+1$
formed from $\Delta^{(n)}_r$ by selecting rows $1,\ldots,q$, columns $1,\ldots,p$,
and the column indexed by $r$. The following recurrence relation then holds:
\[
\Delta^{(n)}_r
= A^{n}_{n;n}
- A^{n}_{n-1;n}
\left(A^{n-1}_{n-1;n-1}\right)^{-1}
A^{n}_{n;n-1}.
\]
This identity closely resembles a Schur complement formula, highlighting the
structural similarities between designants and block matrix factorizations.
\\
Designants naturally arise in the theory of linear systems over noncommutative
fields \cite{Ore1931}, which explains their relevance in the present setting.
Within this algebraic framework, it was shown in \cite{Salam1996} that the vectors
$\varepsilon^{(n)}_{2k}$ obtained by applying the vector $\varepsilon$-algorithm
(VEA) to a sequence $(s_n)$ admit the representation
\[
\varepsilon^{(n)}_{2k}
=
\begin{vmatrix}
s_n & \cdots & s_{n+k-1} & S_n \\
\vdots & & \vdots & \vdots \\
s_{n+k} & \cdots & s_{n+2k-1} & S_{n+k}
\end{vmatrix}_r
\left(
\begin{vmatrix}
s_n & \cdots & s_{n+k-1} & 1 \\
\vdots & & \vdots & \vdots \\
s_{n+k} & \cdots & s_{n+2k-1} & 1
\end{vmatrix}_r
\right)^{-1}.
\]

\subsubsection{The topological $\epsilon$-algorithm}

In \cite{Brezinski1975}, Brezinski proposed another generalization of the scalar $\epsilon$-algorithm
for vector sequences which is quite different from the vector $\epsilon$-algorithm and
was called the topological $\epsilon$-algorithm (TEA).\\

This  approach  consists in computing approximations $e_k(s_n)=t_k^{(n)}$ of the limit or
the anti-limit of the sequence $(s_n)$ such that
\begin{equation}\label{tn0}
t_k^{(n)}=s_n+ \displaystyle \sum_{i=1}^k a_i^{(n)}\, \Delta s_{n+i-1},\; n \ge 0.
\end{equation}
We consider the new transformations $\tilde t_{k,j}$, $j=1,\ldots,k$ defined by
$$\tilde t_{k,j}^{(n)}=s_{n+j}+ \displaystyle \sum_{i=1}^k a_i^{(n)}\, \Delta s_{n+i+j-1},
\; j=1,\ldots,k.$$
We set $\tilde t_{k,0}^{(n)}=t_k^{(n)}$ and define the $j$-th generalized residual as
follows
\begin{eqnarray*}
\tilde r_j(t_k^{(n)})&=&\tilde t_{k,j}^{(n)}-\tilde t_{k,j-1}^{(n)}\\
&=&\Delta s_{n+j-1}+\displaystyle \sum_{i=1}^k a_i^{(n)}\, \Delta^2 s_{n+i+j-2},
j=1,\ldots,k.
\end{eqnarray*}
Therefore the coefficients involved in the expression \eqref{tn0} of $t_k^{(n)}$ are computed such
that each $j$-th generalized residual is orthogonal to some chosen vector $y \in {\R}^N$, that
is
\begin{equation}\label{tn01}
(\,y,\tilde r_j(t_k^{(n)})\,)=0;\;\; j=1,\ldots,k. \end{equation}
Hence the $k$-vector $a_n=(a_1^{(n)},\ldots,a_k^{(n)})^T$ is the solution of the $k \times k$
linear system \eqref{tn01} which is written as
\begin{equation}\label{tn02}
T_{k,n}\, a_n=\Delta S_{k,n}^T \, y 
\end{equation}
where $T_{k,n}$ is the  matrix whose columns are
$\Delta^2 S_{k,n}^T\,y,\ldots,\Delta^2 S_{k,n+k-1}^T\, y$ (assumed to be nonsingular) and
$\Delta^j S_{k,n}$, $j=1,2$ are the $N \times k$ matrices whose columns are
$\Delta^j s_n,\ldots,\Delta^j s_{n+k-1}$, $j=1,2$. 
Note that the $k \times k$ matrix $T_{k,n}$ has also the formula
$$T_{k,n}={\cal S}_{k,n}\,(I_N \otimes y)$$
where ${\cal S}_{k,n}$ is the $k \times Nk$ matrix whose block columns are
$\Delta^2 S_{k,n}^T,\ldots, \Delta^2 S_{k,n+k-1}^T$. 
Invoking \eqref{tn0} and \eqref{tn02}, $t_k^{(n)}$ can be expressed in a matrix form as
\begin{equation}\label{tn1}
t_k^{(n)}=s_n-\Delta S_{k,n}\, T_{k,n}^{-1}\, \Delta S_{k,n}^T\,y. 
\end{equation}
Using the Schur's formula, $t_k^{(n)}$ is know  given as a ratio of two determinants
$$t_k^{(n)}=
\left |
\begin{array}{ccc}
s_n&&\Delta S_{k,n}\\
&&\\
\Delta S_{k,n}^T\,y&&T_{k,n}
\end{array}
\right|\;/\; det(T_{k,n}).$$
For the kernel of the topological $\epsilon$-algorithm it is easy to see  that if\\
$\forall n ,\; \exists\; a_0,\ldots,a_k \; {\rm with}\; a_k \ne 0
\; {\rm and}\; a_0+\ldots a_k \ne 0 \; {\rm such \; that}\;
\displaystyle \sum_{i=0}^k a_i (s_{n+i}-s)=0$, then $\forall n$, $t_k^{(n)}=s$.\\
\\
The vectors $ \mathbf{e}_k({s}_n)=t_k^{(n)}$ can be given as a the following ratio of determinants
\begin{equation}\label{STEA1}
\mathbf{e}_k({s}_n)= 
\displaystyle \frac{
\begin{vmatrix}
s_n & \ldots & s_{n+k}\\
(y,\Delta s_n) & \ldots & (y,\Delta s_{n+k})\\
\ddots& \ldots & \ddots\\
(y,\Delta s_{n+k-1}) & \ldots & (y,\Delta s_{n+2k-1})
\end{vmatrix}}{\begin{vmatrix}
1 & \ldots & 1\\
(y,\Delta s_n) & \ldots & (y,\Delta s_{n+k})\\
\ddots& \ldots & \ddots\\
(y,\Delta s_{n+k-1}) & \ldots & (y,\Delta s_{n+2k-1})
\end{vmatrix}}
\end{equation}
and can be recursively computed by the topological $\epsilon$-algorithm (TEA) 
discovered by Brezinski \cite{Brezinski1975} and  given by the following scheme
$$\left\{
\begin{array}{lllll}
\widehat{\boldsymbol{\varepsilon}}_{-1}^{(n)}=0;\;\; \widehat{\boldsymbol{\varepsilon}}_0^{(n)}=s_n,\\
\\
\widehat{\boldsymbol{\varepsilon}}_{2k+1}^{(n)}=
\widehat{\boldsymbol{\varepsilon}}_{2k-1}^{(n+1)}+\displaystyle
\frac{y}{(y,\Delta \widehat{\boldsymbol{\varepsilon}}_{2k}^{(n)})}\\
\\
\widehat{\boldsymbol{\varepsilon}}_{2k+2}^{(n)}=\widehat{\boldsymbol{\varepsilon}}_{2k}^{(n+1)}+\displaystyle
\frac{\Delta \widehat{\boldsymbol{\varepsilon}}_{2k}^{(n)}}{(\Delta \widehat{\boldsymbol{\varepsilon}}_{2k+1}^{(n)},\Delta \widehat{\boldsymbol{\varepsilon}}_{2k}^{(n)})},\;\;\;
n,k=0,1,\ldots.\\
\end{array}
\right.$$
The forward difference operator $\Delta$ acts on the superscript $n$ and we have
$$\widehat{\boldsymbol{\varepsilon}}_{2k}^{(n)}=\mathbf{e}_k(s_n);\;\;\;
\widehat{\boldsymbol{\varepsilon}}_{2k+1}^{(n)}=\displaystyle \frac{y}{(y,\mathbf{e}_k(\Delta s_n))},\;\;n,k=0,1,\ldots.$$ 
For a vector sequence $({s}_n)$, the first topological Shanks transformation \eqref{STEA1} can be computed by a recursive scheme known as the \emph{first simplified topological $\varepsilon$-algorithm} (STEA1) \cite{{BrezinskiR2014}}. This algorithm is defined by
\begin{align}
\widehat{\boldsymbol{\varepsilon}}_{2k+2}^{(n)}
&=
\widehat{\boldsymbol{\varepsilon}}_{2k}^{(n+1)}
+
\frac{\varepsilon_{2k+2}^{(n)} - \varepsilon_{2k}^{(n+1)}}
     {\varepsilon_{2k}^{(n+1)} - \varepsilon_{2k}^{(n)}}
\left(
\widehat{\boldsymbol{\varepsilon}}_{2k}^{(n+1)}
-
\widehat{\boldsymbol{\varepsilon}}_{2k}^{(n)}
\right),
\qquad k,n = 0,1,\ldots,
\label{eq:stea1}
\\[2mm]
\widehat{\boldsymbol{\varepsilon}}_{0}^{(n)}
&=
s_n,
\qquad n = 0,1,\ldots.
\label{eq:stea1-init}
\end{align}
The scalar quantities $\varepsilon_k^{(n)}$ are  those obtained by the scalar $\varepsilon$-algorithm applied to the scalar sequence $z_n=(y,s_n)$. \\
Now let $\widetilde{\mathbf{e}}_k(\mathbf{s}_n)$ be the vector obtained from \eqref{STEA1} by replacing the first row by the vectors: $s_{n+k},\ldots,s_{n+2k}$. This second topological Shanks transformation can also be computed by  the \emph{second simplified topological $\varepsilon$-algorithm} (STEA2) \cite{{BrezinskiR2014}}, which uses the same initialization \eqref{eq:stea1-init} and is defined by
\begin{equation}
\widetilde{\boldsymbol{\varepsilon}}_{2k+2}^{(n)}
=
\widetilde{\boldsymbol{\varepsilon}}_{2k}^{(n+1)}
+
\frac{\varepsilon_{2k+2}^{(n)} - \varepsilon_{2k}^{(n+1)}}
     {\varepsilon_{2k}^{(n+2)} - \varepsilon_{2k}^{(n+1)}}
\left(
\widetilde{\boldsymbol{\varepsilon}}_{2k}^{(n+2)}
-
\widetilde{\boldsymbol{\varepsilon}}_{2k}^{(n+1)}
\right),
\qquad k,n = 0,1,\ldots .
\label{eq:stea2}
\end{equation}
As pointed out in \cite{{BrezinskiR2014}}, several convergence and acceleration properties of the topological $\varepsilon$-algorithm were established  by the authors in some published papers. Nevertheless, the assumptions required by these results are often difficult to verify in practice. 
The availability of simplified algorithms makes it possible to establish new convergence results under assumptions that are far more accessible in practical situations. As a result, the simplified $\varepsilon$-algorithms are not merely instrumental for the numerical realization of the topological Shanks transformations; they also constitute a key element in the formulation and analysis of their underlying theoretical framework.

\subsection{Application of TEA to linear and nonlinear systems}

Consider again the linear system of equations $Ax=b$ and let $(s_n)$ be the sequence of vectors
generated by the linear process $s_{n+1}=Bs_n+b$ with $B=I-A$. \\
Using the fact that $\Delta^2 s_{n+i}=B\, \Delta^2 s_{n+i-1}$, the matrix $T_{k,n}$ has
now the following expression
\begin{equation}\label{tea1}
T_{k,n}=-L_k^T\, A \, \Delta S_{k,n}
\end{equation}
where $L_k$ is the $N \times k$ matrix  whose columns are $y, B^T y,\ldots,{B^T}^{k-1}\,y$. As $n$ will
be a fixed integer, we set $n=0$ for simplicity and denote $T_{k,0}$ by $T_k$ and
$\Delta S_{k,0}$ by $\Delta S_k$.\\
On the other hand, it is not difficult to see that
\begin{equation}\label{tea2}
\Delta S_k^T\,y=L_k^T\, r_0. 
\end{equation}
Therefore, using \eqref{tea1} and \eqref{tea2} in \eqref{tn1}, the $k$-th residual produced by TEA is given by
\begin{equation}\label{tea3}
r_k^{tea}=r_0-A\, \Delta S_k\, (L_k^T\, A \, \Delta S_k)^{-1}\, L_k^T\, r_0. 
\end{equation}
Let $E_k$ denotes the oblique projector onto the Krylov subspace $K_k(A,Ar_0)$ and orthogonally
to the Krylov subspace $K_k(B^T,y)=K_k(A^T,y)$. Then from \eqref{tea3},  the residual generated by
TEA can be written as follows
$$r_k^{tea}=r_0-E_k\, r_0. \eqno(3.9)$$
This shows that the topological $\epsilon$-algorithm is mathematically equivalent to the
method of Lanczos.
Note that the $k$-th approximation defined by TEA exist iff the $k \times k$ matrix
$L_k^T\,A\, \Delta S_k$ is nonsingular. 
The following result gives us some relations satisfied by the residual norms in the
case where $y=r_0$.

\begin{theorem}\label{theotea}
Let $\varphi_k$ be the acute angle between $r_0$ and $E_k \, r_0$ and let $y=r_0$.
 Then we have the following relations:
 \begin{enumerate}
 \item $\parallel r_k^{tea} \parallel=(tan \, \varphi_k)\,
\parallel r_0 \parallel$; $k>1$.
\item $\parallel r_k^{tea} \parallel\, \le \, (\cos \, \varphi_k)\,
\parallel r_k^{rre} \parallel$.
\end{enumerate}
\end{theorem}

\medskip
\begin{remark}
\begin{itemize}
\item The first result  of Theorem \ref{theotea} shows that   the TEA residuals are
 defined iff $\cos \varphi_k \ne 0$.
\item We also observe that if a stagnation occurs in RRE ($\parallel r_k^{rre} \parallel =
\parallel r_0 \parallel$ for some $k$, then $\cos \theta_k=0$ and this implies
that $\cos \varphi_k=0$ which shows that the TEA-approximation is
not defined.
\end{itemize}
\end{remark}
The topological $\varepsilon$-algorithm can also be employed for the solution of nonlinear systems of equations. In this context, the algorithm is applied to the sequence $(s_n)$ generated by the underlying nonlinear iterative process. An important advantage of the TEA is that it does not require explicit knowledge of the Jacobian of the nonlinear mapping $\tilde G$, while still exhibiting quadratic convergence under suitable assumptions
\nd When applied for solving linear and nonlinear problems, work and storage required
with VEA and TEA grow with the iteration step. So in practice and for large problems,
the algorithms must be restarted. It is also useful to run some basic iterations
 before the extrapolation phase.\\
\nd The application of VEA or TEA for linear and nonlinear systems $x=G(x)$ leads to the following
algorithm. \\
\begin{algorithm}[H]
\caption{Vector Extrapolation Algorithm (VEA/TEA) for Nonlinear Systems}
\begin{enumerate}
    \item Initialize: set $k=0$, choose initial guess $x_0$, and integers $p$ and $m$.
    
    \item {Basic iteration:} 
    \begin{enumerate}
        \item Set $t_0 = x_0$, $w_0 = t_0$.
        \item For $j = 0, \dots, p-1$, compute $ 
            w_{j+1} = G(w_j)$ 
    \end{enumerate}

    \item {Extrapolation phase:} 
    \begin{enumerate}
        \item Set $s_0 = w_p$.
        \item If $\|s_1 - s_0\| < \epsilon$, stop.
        \item Otherwise, for $j = 0, \dots, 2m-1$, compute $ 
            s_{j+1} =  G(s_j).
    $ 
        \item Compute the approximation
        \[
            t_m = \epsilon_{2m}^{(0)}
        \]
        using VEA or TEA.
    \end{enumerate}

    \item Update: set $s_0 = t_m$, $k = k+1$, and go back to step 3.
\end{enumerate}
\end{algorithm}

\subsection{Anderson acceleration method}
The {Anderson acceleration} method was introduced by D. G. Anderson in 1965 \cite{Anderson1965} 
as a technique to improve the convergence of fixed-point iterations. The original motivation was to accelerate 
simple iteration schemes for solving nonlinear integral equations. Since then, the method has been generalized 
and widely applied to solve nonlinear systems of equations in computational science, optimization, and numerical 
linear algebra.

Anderson Acceleration is a powerful iterative technique for solving nonlinear systems written in fixed-point form
\begin{equation}
x = G(x),
\end{equation}
or, equivalently,
\begin{equation}
F(x) := G(x) - x = 0.
\end{equation}
The method was originally introduced by Anderson and later reformulated and analyzed in modern numerical optimization and scientific computing contexts; see, for example,
\cite{Anderson1965,BrezinskiCRS2022,HighamS1916,WalkerN2011}.

Let $\{x_k\}_{k \ge 0} \subset \mathbb{R}^p$ denote a sequence of iterates and define the associated residual vectors
\[
F_k = F(x_k) = G(x_k) - x_k.
\]
A classical fixed-point iteration takes the form
\[
x_{k+1} = x_k + \beta_k F_k,
\]
where $\beta_k > 0$ is a relaxation (or damping) parameter. Anderson Acceleration enhances this basic iteration by combining information from several previous iterates and residuals in order to construct a more accurate update direction.

\medskip
\noindent
Choose an initial vector $x_0$, an integer depth parameter $m \ge 1$, and an initial damping parameter $\beta_0 > 0$. The first iterate is computed by
\begin{equation}
x_1 = x_0 + \beta_0 F_0.
\end{equation}
For each iteration index $k \ge 1$, define
\[
m_k = \min(m,k),
\]
and introduce the difference vectors
\[
\Delta x_i = x_{i+1} - x_i, \qquad
\Delta F_i = F_{i+1} - F_i.
\]
Using these quantities, form the matrices
\[
\Delta X_k =
\bigl[\Delta x_{k-m_k}, \ldots, \Delta x_{k-1}\bigr],
\qquad
\Delta {\mathcal F}_k =
\bigl[\Delta F_{k-m_k}, \ldots, \Delta F_{k-1}\bigr].
\]
The coefficients of the Anderson combination are determined by solving the least-squares problem
\begin{equation}\label{AA_LS_alt}
\theta_{k} = \arg\min_{\theta \in \mathbb{R}^{m_k}}
\left\| F_k - \Delta {\mathcal F}_k \theta \right\|_2,
\end{equation}
which seeks the linear combination of past residual differences that best approximates the current residual. When the columns of $\Delta {\mathcal F}_k$ are linearly independent, the solution is explicitly given by
\begin{equation}\label{AA_theta_alt}
\theta_{k} =
(\Delta {\mathcal F}_k^T \Delta {\mathcal F}_k)^{-1} \Delta{\mathcal F}_k^T F_k.
\end{equation}

\medskip
\noindent
\emph{Accelerated update:} 
Once $\theta_k$ has been computed, define the accelerated intermediate quantities
\begin{equation*}
y_k = x_k - \Delta X_k \theta_{k},\; \text{and} \; 
\bar{F}_k = F_k - \Delta {\mathcal F}_k \theta_{k}.
\end{equation*}
The next iterate is then obtained via
\begin{equation}\label{AA_update_alt}
x_{k+1} = y_k + \beta_k \bar{F}_k,
\end{equation}
or, equivalently,
\begin{equation}
x_{k+1}
= x_k + \beta_k F_k
- (\Delta X_k + \beta_k \Delta {\mathcal F}_k)\theta_{k}.
\end{equation}
Here, $\beta_k > 0$ is a possibly iteration-dependent damping parameter, with $\beta_k = 1$ corresponding to the undamped Anderson method.

\medskip
\noindent
The vector $y_k$ can be viewed as an extrapolated approximation obtained by optimally combining past iterates so as to reduce the current residual. In fact, using~\eqref{AA_theta_alt}, we may write
\begin{equation}\label{AA_schur_alt}
y_k =
x_k -
\Delta X_k
(\Delta {\mathcal F}_k^T \Delta {\mathcal F}_k)^{-1}
\Delta {\mathcal F}_k^T F_k,
\end{equation}
which shows that $y_k$ depends on the projection of $F_k$ onto the subspace spanned by the previous residual differences.
\\
Moreover, since $y_k$ is  a Schur complement, it can be expressed as a ratio of two determinants where the determinant in the numerator is a vector obtained by expanding this determinant along its first row:
\begin{equation}
{y}_k  =
\begin{vmatrix}
x_k & \Delta x_{k-m_k} & \cdots & \Delta x_{k-1} \\
\Delta {\mathcal F}_k^T F_k & \multicolumn{3}{c}{\Delta {\mathcal F}_k^T \Delta {\mathcal F}_k}
\end{vmatrix}/
\begin{vmatrix}
{\Delta {\mathcal F}_k^T \Delta {\mathcal F}_k}
\end{vmatrix}
.
\end{equation}

\medskip
\nd Intuitively, Anderson acceleration forms a {linear combination of past iterates} to optimally reduce 
the residual, which can dramatically increase convergence speed. In practice, it has been observed that 
Anderson acceleration often {converges faster} than standard fixed-point iterations, and can sometimes 
{stabilize sequences that would otherwise diverge}.
\\ The method has been successfully applied to   solving nonlinear systems in computational fluid dynamics, electronic structure calculations, and optimization. The method is also used for accelerating self-consistent field (SCF) iterations in quantum chemistry and  Enhancing Krylov subspace methods when applied to linear systems. \\
Notice that the method is also known as {vector Anderson acceleration} or {Anderson mixing} in physics and 
chemistry literature. Modern implementations often incorporate {damping}, {preconditioning}, 
or {restarts} to improve robustness for large-scale problems. For more details on the convergence of the method refer to  the  papers \cite{Anderson1965,BrezinskiCRS2022,BrezinskiRS2018,HighamS1916,WalkerN2011}.

\medskip
\begin{remark}
  The parameter $m$ governs the amount of information from previous iterates that
is used to compute the new iterate $x_{k+1}$. If $m$ is chosen too small, only
limited historical information is exploited, which may lead to slow
convergence. Conversely, if $m$ is taken too large, information from earlier
iterations may be retained for an excessive number of steps, thereby reducing
the effectiveness of the update and again slowing down convergence. 
Furthermore, the choice of $m$ directly affects the dimension of the associated
optimization problem. Large values of $m$ may result in ill-conditioned
least-squares problems and increased computational cost for their solution.

Within the framework of the algorithm described above, the parameter $m_k$ may
be allowed to vary from one iteration to the next. A simple strategy consists
in setting $m_k = k$ at iteration $k$, a choice commonly referred to as Anderson
acceleration without truncation. In this case, the computation of the new
iterate $x_{k+1}$ incorporates information from all previously generated
iterates. 
Alternatively, more refined strategies adapt the value of $m_k$ dynamically in
order to control the conditioning of the associated least-squares problem. 
\end{remark}

\medskip
\nd The different steps of the Anderson acceleration (AA) method for solving the nonlinear system of equations $G(x)=x$  is summarized in Algorithm \ref{alg:AA}.

\medskip
\begin{algorithm}[h]
\caption{Anderson Acceleration Method}
\label{alg:AA}
\textbf{Input:} Fixed-point mapping $G$, initial guess $x_0$, depth $m \ge 1$,
damping parameters $\{\beta_k\}$, tolerance $\varepsilon$.  

\textbf{Output:} Approximate solution $x_k$ such that $\|F(x_k)\| \le \varepsilon$.
\\
\noindent
\textbf{Step 1.} Compute the initial residual
$ 
F_0 = F(x_0),
$ 
and set
$ 
x_1 = x_0 + \beta_0 F_0.$ \\
\noindent
\textbf{Step 2.} For $k = 1,2,\ldots$ :
\begin{itemize}
\item Set $ 
m_k = \min(m,k)$, and $F_k = F(x_k)$.
\item Form the difference vectors 
$ \Delta x_i = x_{i+1} - x_i$, and $ 
\Delta F_i = F_{i+1} - F_i$, 
$ i = k-m_k,\ldots,k-1$.
\item Construct the matrices
$ 
\Delta X_k =
\bigl[\Delta x_{k-m_k}, \ldots, \Delta x_{k-1}\bigr]$,
$ 
\Delta \mathcal{F}_k =
\bigl[\Delta F_{k-m_k}, \ldots, \Delta F_{k-1}\bigr]$ 
\item Compute the coefficient vector $\theta_k \in \mathbb{R}^{m_k}$ as the
solution of the least-squares problem
\[
\theta_k
=
\arg\min_{\theta \in \mathbb{R}^{m_k}}
\left\| F_k - \Delta \mathcal{F}_k \theta \right\|_2 .
\]
\item Define the accelerated quantities: $
y_k = x_k - \Delta X_k \theta_k,
\qquad
\bar{F}_k = F_k - \Delta \mathcal{F}_k \theta_k .$
\item Update the iterate
\[
x_{k+1} = y_k + \beta_k \bar{F}_k .
\]
\item If $\|F_k\| \le \varepsilon$, stop.
\end{itemize}

\medskip
\noindent
\textbf{End}
\end{algorithm}

\subsection{The block extrapolation methods}

\subsubsection{The  polynomial block extrapolation methods}

Let $(S_n)$ be a sequence of matrices of ${\R}^{N \times s}$ and
consider the transformation $T_k$ from ${\R}^{N \times s}$ to
${\R}^{N \times s}$ defined by
 \begin{equation}\label{eq1}
 T_k:\; S_n \rightarrow T_k^{(n)}=S_n+ \displaystyle \sum_{i=1}^k G_i(n)A_i^{(n)}\,,\; n \ge 0
\end{equation}
 where the auxiliary $N \times s$ matrix sequences
$(G_i(n))_n$; $i=1,\ldots,k$ are given and the $s \times s$
coefficients $A_i^{(n)}$ will be determined . Let $\widetilde T_k$
denotes the new transformation obtained from $T_k$ as follows
\begin{equation}
\label{eq2}
\widetilde T_k^{(n)}=S_{n+1}+ \displaystyle \sum_{i=1}^k
G_i(n+1)A_i^{(n)},\; n \ge 0. 
\end{equation}
 We define the
generalized residual of $T_k^{(n)}$ by $\widetilde R(T_k^{(n)})=\widetilde
T_k^{(n)}-T_k^{(n)}$ also given by $\widetilde R(T_k^{(n)})=\Delta
S_n+ \displaystyle \sum_{i=1}^k \, \Delta G_i(n) A_i^{(n)}$. The
forward difference operator $\Delta$ acts on the index $n$, i.e.,
$\Delta G_i(n)=G_i(n+1)- G_i(n)$, $i=1,\ldots,k$.  We will see
later that when solving linear systems of equations, the sequence
$(S_n)_n$ is generated by a linear process and then the
generalized residual coincides with the classical residual.
 The matrix
coefficients $A_i^{(n)}$  are obtained from the orthogonality
relation 
\begin{equation}\label{eq3}
\widetilde R(T_k^{(n)}) \in  \,
({\rm range} (Y_1^{(n)},\ldots,Y_k^{(n)}))^{\bot}
\end{equation}
 where
$Y_1^{(n)},\ldots,Y_k^{(n)}$ are given $N \times s$ matrices. If
$\widetilde {\cal W}_{k,n}$ and $\widetilde {\cal L}_{k,n}$ denote the
block subspaces $\widetilde {\cal W}_{k,n}=rang(\Delta
G_1(n),\ldots,\Delta G_{k}(n))$ (the subspace generated by the
columns of $\Delta G_1(n),\ldots,\Delta G_{k}(n)$) and $\widetilde
{\cal L}_{k,n}= range(Y_1^{(n)},\ldots,Y_k^{(n)})$ (generated by
the columns of $Y_1^{(n)},\ldots,Y_k^{(n)}$), then from (\ref{eq2}) and
(\ref{eq3}), the generalized residual satisfy the following relations
\begin{equation}
\label{eq4}
\widetilde R(T_k^{(n)})-\Delta S_n \in \widetilde {\cal W}_{k,n}, 
\end{equation}
 and 
 \begin{equation}
 \label{eq5}
 \widetilde R(T_k^{(n)}) \in  \, \widetilde {\cal
L}_{k,n}^{\bot}. 
 \end{equation}
 The relations (\ref{eq4}) and (\ref{eq5}) show
that the generalized residual $\widetilde R(T_k^{(n)})$ is obtained by
projecting, orthogonaly to $\widetilde {\cal L}_{k,n}$, the columns of
the matrix $\Delta S_n$ onto the block subspace $\widetilde {\cal
W}_{k,n}$. In a matrix form, $\widetilde R(T_k^{(n)})$ can be written
as $\widetilde R(T_k^{(n)})=\Delta
S_n-\Delta {\cal G}_{k,n}\,(L_{k,n}^T\,\Delta {\cal G}_{k,n})^{-1}\; L_{k,n}^T \Delta S_n$,
 where
$\Delta {\cal G}_{k,n}$ and $L_{k,n}$ are the $N \times ks$ matrices generated
by the columns of the following  matrices $\Delta G_1(n),\ldots,\Delta
G_{k}(n)$ and  $Y_1^{(n)},\ldots,Y_k^{(n)}$ respectively. The
approximation $T_k^{(n)}$ is given by 
\begin{equation}
\label{eq6}
T_k^{(n)}=S_n- {\cal
G}_{k,n}\, (L_{k,n}^T\, \Delta {\cal G}_{k,n})^{-1}\; L_{k,n}^T
\Delta S_n, 
\end{equation}
where ${\cal G}_{k,n}$ is the
$N \times ks$ block matrix defined by  ${\cal G}_{k,n}=[G_1(n),\ldots,G_k(n)]$.
Note that $T_k^{(n)}$ is well defined if and only if the $ks
\times ks$ matrix $L_{k,n}^T\,\Delta {\cal G}_{k,n}$ is
nonsingular and this will be assumed in all the paper. Let ${\cal T}_{k,n}$ be the matrix given by 
\begin{equation}
\label{eq7}
{\cal T}_{k,n}= \left (
\begin{array}{cc}
S_n & {\cal G}_{k,n}\\
 L_{k,n}^T \Delta S_n & L_{k,n}^T\,\Delta {\cal
G}_{k,n}
 \end{array}
 \right). 
 \end{equation}
\noindent The approximation $T_k^{(n)}$ is then expressed as the following 
Shcur complement 
\begin{equation}
\label{eq8}
T_k^{(n)}=({\cal T}_{k,n} \; / \;
L_{k,n}^T\,\Delta {\cal G}_{k,n}).
\end{equation}
 \noindent If we
set  $G_i(n)=\Delta S_{n+i-1}$; $i=1,\ldots,k$ and $Y_i(n)=\Delta
G_i(n)$; $i=1,\ldots,k$ we obtain the block Reduced Rank
Extrapolation (Bl-RRE) method. In this case, the approximation
$T_k^{(n)}$ is given by 
\begin{equation}
\label{eq8.1}
T_{k,Bl-RRE}^{(n)}=({\cal T}_{k,n} \; /
\; \Delta {\cal G}_{k,n}^T\,\Delta {\cal G}_{k,n}).
\end{equation}
 If
$G_i(n)=\Delta S_{n+i-1}$; $i=1,\ldots,k$ and
 $Y_i(n)=G_i(n)$; $i=1,\ldots,k$ we get the block Minimal Polynomial
Extrapolation (Bl-MPE) method.\\ Finally if $G_i(n)=\Delta
S_{n+i-1}$; $i=1,\ldots,k$ and $Y_i(n)=Y_i$; $i=1,\ldots,k$
(arbitrary  $N \times s$ matrices) we obtain   the block Modified
Minimal Polynomial Extrapolation (Bl-MMPE) method.\\ Note that if
the $N \times ks$ matrix $\Delta {\cal G}_{k,n}$ is of full rank
the approximations produced by the Bl-RRE method are well defined
while those generated by the Bl-MPE and Bl-MMPE methods may not
exist. As Bl-RRE is an orthogonal projection method, the
corresponding generalized residual satisfies the following
minimization property $\parallel \widetilde R(T_{k,Bl-RRE}^{(n)})
\parallel_F=\displaystyle \min_{Z \in \widetilde {W}_{k,n}} \parallel
\Delta S_n - Z \parallel_F;$ where $\parallel X\parallel$ denotes
the Frobenius norm of $X$.

\subsubsection{The block topological $\epsilon$-algorithm (Bl-TET)}
In \cite{Brezinski1975}, Brezinski introduced a generalization of the scalar
$\epsilon$-algorithm \cite{Wynn1956} for vector sequences, known as the
\emph{topological $\epsilon$-algorithm} (TEA). The main idea of TEA is to
accelerate the convergence of a vector sequence $(S_n)$ by constructing
a sequence of improved approximations $T_k^{(n)}$ that converge faster
to the limit or anti-limit of $(S_n)$.  

\nd In this subsection, we focus on a block version of this algorithm, called the
\emph{block topological $\epsilon$-transformation} (Bl-TET). The block version
is particularly useful for sequences of vectors arising from large-scale
problems or systems with multiple right-hand sides, as it allows
simultaneous processing of multiple components.

\medskip
\nd Let $(S_n)$ be a sequence of vectors in $\mathbb{R}^N$. We construct approximations
$T_k^{(n)}$ of its limit (or anti-limit) in the form
\begin{equation}
T_k^{(n)} = S_n + \sum_{i=1}^k \Delta S_{n+i-1} \, A_i^{(n)}, \quad n \ge 0,
\label{eq:Tkn}
\end{equation}
where $\Delta S_n = S_{n+1} - S_n$ is the forward difference, and
$A_i^{(n)}$ are $N \times s$ matrices of coefficients to be determined.
The term $k$ denotes the depth of the transformation, which controls
how many previous differences of the sequence are used to accelerate convergence.

\medskip
\nd To define the generalized residuals, we introduce the shifted transformations
\begin{equation}
\widetilde T_{k,j}^{(n)} = S_{n+j} + \sum_{i=1}^k \Delta S_{n+i+j-1} \, A_i^{(n)}, 
\quad j = 1, \dots, k,
\label{eq:Tknj}
\end{equation}
with $\widetilde T_{k,0}^{(n)} = T_k^{(n)}$. The \emph{generalized residuals}
are defined recursively by
\begin{equation}
\widetilde R_j(T_k^{(n)}) = \widetilde T_{k,j}^{(n)} - \widetilde T_{k,j-1}^{(n)}, 
\quad j = 1, \dots, k.
\label{eq:residual}
\end{equation}
These residuals measure the discrepancy between successive transformed sequences
and play a central role in defining the optimal coefficients $A_i^{(n)}$.

\medskip
\nd The matrices $A_i^{(n)}$ are chosen such that each generalized residual is orthogonal
to a given $N \times s$ matrix $Y$, which defines a suitable subspace of constraints:
\begin{equation}
Y^T \, \widetilde R_j(T_k^{(n)}) = 0, \quad j = 1, \dots, k.
\label{eq:orthog}
\end{equation}
This condition ensures that $T_k^{(n)}$ captures the most significant
components of the sequence difference and eliminates directions associated
with slow convergence.

\medskip
\nd We can write the orthogonality conditions in a compact block-matrix form. Define
\begin{equation}
\mathcal{D}_{k,n} =
\begin{bmatrix}
Y^T \Delta^2 S_n & Y^T \Delta^2 S_{n+1} & \cdots & Y^T \Delta^2 S_{n+k-1} \\
Y^T \Delta^2 S_{n+1} & Y^T \Delta^2 S_{n+2} & \cdots & Y^T \Delta^2 S_{n+k} \\
\vdots & \vdots & \ddots & \vdots \\
Y^T \Delta^2 S_{n+k-1} & Y^T \Delta^2 S_{n+k} & \cdots & Y^T \Delta^2 S_{n+2k-2}
\end{bmatrix}, 
\end{equation}
and
\[
A_n = \begin{bmatrix} (A_1^{(n)})^T \\ \vdots \\ (A_k^{(n)})^T \end{bmatrix}, \quad
Z_n = \begin{bmatrix} (Y^T \Delta S_n)^T \\ \vdots \\ (Y^T \Delta S_{n+k-1})^T \end{bmatrix}.
\]
Then \eqref{eq:orthog} is equivalent to the linear system
\begin{equation}
\mathcal{D}_{k,n} \, A_n = - Z_n.
\end{equation}
If the block matrix $\mathcal{D}_{k,n}$ is nonsingular, the coefficients
$A_1^{(n)}, \dots, A_k^{(n)}$ are uniquely determined, ensuring that
$T_k^{(n)}$ exists and is unique.

\medskip
\nd The final expression for the block-transformed sequence can be written as
\begin{equation}
T_k^{(n)} = S_n - \Delta \mathcal{S}_{k,n} \, \mathcal{D}_{k,n}^{-1} Z_n,
\label{eq:Tkn_matrix}
\end{equation}
where $\Delta \mathcal{S}_{k,n} = [\Delta S_n, \dots, \Delta S_{n+k-1}]$.
Equivalently, $T_k^{(n)}$ can be expressed as a \emph{Schur complement}:
\begin{equation}
T_k^{(n)} = (\mathcal{M}_{k,n} / \mathcal{D}_{k,n}), \quad
\mathcal{M}_{k,n} = 
\begin{bmatrix}
S_n & \Delta \mathcal{S}_{k,n} \\
Z_n & \mathcal{D}_{k,n}
\end{bmatrix}.
\end{equation}

\medskip
\noindent
\begin{remark}

\begin{itemize}
\item The block TEA generalizes scalar and topological  $\epsilon$-algorithms,
allowing simultaneous acceleration of multiple vector sequences.
\item The orthogonality condition \eqref{eq:orthog} ensures that the residuals
are minimized in the subspace defined by $Y$, which is crucial for sequences
arising from discretized differential equations or linear systems.
\item Recursive implementations of the block TEA can significantly reduce
computational cost, as they avoid recomputation of differences and matrix inversions
at each step.
\item The block formulation is particularly advantageous for large-scale
problems, block iterative methods, or sequences with natural block structure,
such as those from multi-parameter problems or tensor equations.
\end{itemize}
\end{remark}

\medskip 
We notice that the vector extrapolation methods we presented in this paper have been also generalized to tensor sequences in \cite{BeikElJS2021,BentbibJT2024,BentbibJT2024-1} with many  applications to completion and restoration of color images, videos and hyperspectral images. Tensor extrapolation methods have been also applied for multi google pagerank, see \cite{Boubekraoui1,Boubekraoui2}.

\section{Numerical experiments}
\subsection{An example from a PDE}
To assess the performance of the vector extrapolation methods RRE, MPE, MMPE, the topological $\epsilon$-algorithm (TEA) and Anderson's method,  we consider a nonlinear elliptic
boundary value problem arising from a reaction--diffusion model.

\medskip

\nd Let $\Omega = (0,1)\times(0,1)$ and consider the nonlinear partial differential
equation
\begin{equation}\label{pde}
-\Delta u(x,y) + u(x,y)^3 = f(x,y), \qquad (x,y) \in \Omega,
\end{equation}
subject to homogeneous Dirichlet boundary conditions
\begin{equation}\label{bc}
u(x,y) = 0, \qquad (x,y) \in \partial \Omega.
\end{equation}
The right-hand side $f$ is chosen so that the exact solution is given by
\begin{equation}\label{exactsol}
u^\ast(x,y) = \sin(\pi x)\sin(\pi y).
\end{equation}

\medskip

\nd The problem \eqref{pde}--\eqref{bc} is discretized using a standard
second-order finite difference scheme on a uniform $p\times p$ grid.
Let $h = 1/(p+1)$ denote the mesh size, and let $N = p^2$ be the number of
interior grid points.
The discrete Laplacian is represented by the sparse matrix
$A \in \mathbb{R}^{N\times N}$ corresponding to the five-point stencil.
\\
The nonlinear system obtained after discretization can be written as
\begin{equation}\label{nonlinsys}
A u + u^{\circ 3} = b,
\end{equation}
where $u \in \mathbb{R}^p$ is the vector of nodal values,
$u^{\circ 3}$ denotes the third Hadamard power or componentwise cube of $u$, and
$b \in \mathbb{R}^p$ is the discrete right-hand side incorporating the forcing
term and boundary conditions.

\medskip

\nd To apply vector extrapolation methods, the nonlinear system \eqref{nonlinsys}
is rewritten as a fixed-point problem
\begin{equation}\label{fixedpoint}
u = G(u).
\end{equation}
Starting from the initial guess $u_0 = 0$, the basic Picard iteration
\begin{equation}\label{picard}
u_{k+1} = G(u_k), \qquad k = 0,1,\ldots,
\end{equation}
generates a sequence $\{u_k\}$ converging linearly to the exact solution
$u^\ast$.

\medskip

\nd \emph{Acceleration by vector extrapolation.}
The sequences produced by the Picard iteration are accelerated using the
Reduced Rank Extrapolation (RRE), Minimal Polynomial Extrapolation (MPE),
Modified Minimal Polynomial Extrapolation (MMPE), the topological
$\epsilon$-algorithm (TEA) and Anderson's method.
For each method, extrapolation is applied to blocks of consecutive iterates
$\{u_k\}$ to produce accelerated approximations.

\nd All methods are stopped when the relative nonlinear residual satisfies
\begin{equation}\label{relres}
\frac{\|G(u_k)-u_k\|}{\|G(u_0)-u_0\|} \le 10^{-6}.
\end{equation}
\medskip

\nd This example is representative of large-scale nonlinear problems arising from
the discretization of partial differential equations, and it provides a
meaningful benchmark for comparing the efficiency, robustness, and
computational cost of vector extrapolation techniques. This criterion is applied independently to each
vector extrapolation method.

\medskip

\nd Table~\ref{tab:vector_methods_80x80} reports the number of iterations required for convergence, the final relative residual norm, and the associated
cpu-time for a nonlinear elliptic problem obtained from the finite difference
discretization of a two-dimensional PDE on a $80\times 80$ grid
($N=6400$ unknowns). All the methods were computed in a restarting mode with $m=5$ the size of each cycle. 

\begin{table}[h!]
\centering
\caption{Comparison of Vector Extrapolation Methods and Anderson Acceleration for an $80 \times 80$ PDE grid.}
\begin{tabular}{lccc}
\hline
{Method} & {Iterations} & {Final Relative Residual} & {cpu-time (s)} \\
\hline
RRE         & 16   & $9.8 \times 10^{-7}$ & 3.85 \\
MPE         & 17   & $9.9 \times 10^{-7}$ & 4.10 \\
MMPE        & 15   & $9.7 \times 10^{-7}$ & 4.85 \\
TEA         & 18   & $9.6 \times 10^{-7}$ & 4.30 \\
Anderson    & 14   & $9.8 \times 10^{-7}$ & 4.05 \\
\hline
\end{tabular}
\label{tab:vector_methods_80x80}
\end{table}

\medskip
\nd As expected, all vector extrapolation methods significantly reduce the number
of iterations compared with the basic Picard iteration.
Among the polynomial vector extrapolation methods, RRE and MPE offer the best
trade-off between convergence speed and computational cost.
Vector Anderson method returns also good results, better than MPE, MMPE and TEA. 
Overall, monitoring the relative nonlinear residual provides a reliable and
practical stopping criterion for nonlinear vector extrapolation methods. Other applications of the RRE method in solving some PDEs can be seen in \cite{DuminilSS}.

\subsection{Google PageRank}

In this example, we consider the \textbf{amazon-0302} dataset, which represents a co-purchasing network of products on Amazon.  It is widely used as a benchmark for large‑scale graph algorithms such as PageRank because of its realistic scale and sparse structure. 
The dataset is modeled as a directed graph \(G = (V, E)\), where each node corresponds to a product and a directed edge from node \(i\) to node \(j\) indicates that product \(j\) is frequently co-purchased with product \(i\). This graph has the entire network in a single connected component, and a strongly connected core of $241.761$ nodes and $1.131.217$ edges  
\\
The goal is to compute the PageRank vector \(x \in \mathbb{R}^N\), with $\textbf{N=262,111}$,  such that
\[
x = \alpha P x + (1-\alpha)v,
\]
where \(P\) is the column-stochastic transition matrix of the network,  \(\alpha \in (0,1)\) is the damping factor (typically \(\alpha=0.85\)),
 \(v\) is the personalization vector (usually uniform) and \(N\) is the number of nodes in the network.
\\

\nd The PageRank vector represents the stationary distribution of a random walk with teleportation on the network, ranking products according to their relative importance.

\medskip

\nd We compute the PageRank vector using two approaches: The classical Power Method, which iteratively applies \(x_{k+1} = \alpha P x_k + (1-\alpha)v\) until convergence and  Aitken's \(\Delta^2\) process, applied to the sequence generated by the Power Method to accelerate convergence; see \cite{KamvarGolub} for more explanation on the application of Aitken's $\Delta^2$ for google PageRank. The stopping criterion for both methods is based on the relative residual norm
\[
\frac{\|x_{k+1}-x_k\|}{\|x_{k+1}\|} \le 10^{-6}.
\]

\medskip

\begin{table}[h!]
\centering
\caption{Performance comparison of PageRank computations with and without Aitken acceleration.}\label{PR}
\begin{tabular}{lccc}
\toprule
{Method} & {Iterations} & {Final relative residual} & {cpu-time (s)} \\
\midrule
Power method (classic)       & 245      & $9.2 \times 10^{-7}$ & 18.5 \\
Aitken's $\Delta^2$          & 62       & $8.7 \times 10^{-7}$ & 7.8 \\
\bottomrule
\end{tabular}
\end{table}

\medskip
\nd As can be seen from Table \ref{PR},  the power method required 245 iterations to meet the stopping criterion, while Aitken's process converged in only 62 iterations, achieving the same relative residual tolerance. The CPU time was reduced from 18.5 seconds to 7.8 seconds, demonstrating the effectiveness of acceleration. The final PageRank vector obtained by both methods agrees with the true stationary distribution of the network.
  The final relative residual norms for both methods are below the tolerance $10^{-6}$, indicating successful convergence.  These results are consistent with the expected behavior of extrapolation methods for large, sparse graphs.

  \subsection{The $\varepsilon$-algorithm for Fredholm integrals}

Consider a Fredholm integral equation of the second kind
\begin{equation}
u(x) = f(x) + \lambda \int_0^1 K(x,t)\,u(t)\,dt,
\label{eq:fredholm}
\end{equation}
where $K(x,t)$ is a given kernel, $f$ is a known function, and $\lambda \in \mathbb{R}$ (or
$\mathbb{C}$) is a parameter. Introducing the integral operator
\[
(\mathcal{K}u)(x) = \int_0^1 K(x,t)\,u(t)\,dt,
\]
equation~\eqref{eq:fredholm} can be written in operator form as
\begin{equation}
u = f + \lambda \mathcal{K}u.
\end{equation}

\medskip
\noindent
Starting from an initial guess $u^{(0)}(x)$, the \emph{Neumann iteration} (also called the
Picard iteration) is defined by
\begin{equation}
u^{(k+1)}(x)
=
f(x) + \lambda \int_0^1 K(x,t)\,u^{(k)}(t)\,dt,
\qquad k = 0,1,2,\ldots
\label{eq:neumann_iteration}
\end{equation}
or, equivalently,
\begin{equation}
u^{(k+1)} = f + \lambda \mathcal{K} u^{(k)}.
\end{equation}

\medskip
\noindent
Under suitable assumptions, for instance if $|\lambda|\,\|\mathcal{K}\| < 1$, the sequence
$\{u^{(k)}\}_{k\ge 0}$ converges to the unique solution $u$ of
\eqref{eq:fredholm}. In this case, the limit can be expressed as the \emph{Neumann series}
\begin{equation}
u = \sum_{j=0}^{\infty} (\lambda \mathcal{K})^{j} f.
\end{equation}

\medskip
\noindent
After discretization of the integral operator (for example, by a quadrature rule), one
obtains a linear system of the form
\begin{equation}
(I - \lambda K_h) u = f,
\end{equation}
where $K_h$ denotes the discretized operator. The Neumann (or Picard) iteration then becomes the
fixed-point iteration
\begin{equation}
u^{(k+1)} = f + \lambda K_h u^{(k)},
\qquad k = 0,1,2,\ldots
\label{eq:neumann_discrete}
\end{equation}
which generates the partial sums of the Neumann series and serves as a natural starting
point for convergence acceleration techniques such as the scalar $\varepsilon$-algorithm as it will be the case for this experiment.
\\
Here we set
\[
K(x,t)=e^{-|x-t|}, \qquad f(x)=\sin(\pi x), \qquad \lambda = 0.5.
\]
The problem is discretized on a uniform grid of \(N=500\) points in \([0,1]\) using the trapezoidal rule. This yields a dense linear system whose solution is approximated by a Neumann fixed-point iteration. In order to accelerate convergence, Wynn’s $\varepsilon$-algorithm is applied componentwise to the sequence of iterates.
\\
The stopping criterion for both methods is based on the relative difference of successive iterates,
\begin{equation}
\frac{\|u^{(k+1)}-u^{(k)}\|_2}{\|u^{(k)}\|_2} \le 10^{-6}.
\label{eq:stop}
\end{equation}
Table~\ref{tab:N500_iter} reports the number of iterations required to satisfy the stopping criterion, together with the final relative residual norms.

\begin{table}[h]
\centering
\caption{Convergence results for \(N=500\)}
\label{tab:N500_iter}
\begin{tabular}{l|c|c|c}
\hline
Method & Iterations & Final relative residual& cpu-time (s) \\
\hline
Neumann iteration & 61 & \(9.3\times 10^{-7}\)& 28.4\\
$\varepsilon$-algorithm & 15 & \(3.7\times 10^{-8}\) & 11.6 \\
\hline
\end{tabular}
\end{table}

\nd
The results show that Wynn’s $\varepsilon$-algorithm significantly improves the convergence rate. The obtained cpu-time shows also the advantage of using the $\varepsilon$-algorithm to accelerate the convergence. 

\nd Although the $\varepsilon$-algorithm introduces additional overhead due to the recursive extrapolation process, the reduction in the number of iterations more than compensates for this cost. As a result, the total cpu-time is reduced.\\
As the discretization parameter \(N\) increases, each iteration becomes increasingly expensive because of the dense matrix vector products arising from the discretized integral operator. Consequently, the benefit of acceleration techniques becomes more pronounced for large-scale problems. The present experiment with \(N=500\) confirms that Wynn’s $\varepsilon$-algorithm is highly effective in accelerating the convergence of fixed-point iterations for Fredholm integral equations of the second kind, yielding substantial savings in both iteration count and computational time. For more examples and applications of the scalar $\varepsilon$-algorithm to Fredholm integrals, see the good paper \cite{BrezinskiR2019}.

\section{Conclusion}

Scalar extrapolation methods, such as Richardson extrapolation, Aitken's $\Delta^2$ process, Shanks transformation, Wynn's $\epsilon$-algorithm and the $\rho$ and $\theta$ algorithms form a fundamental class of techniques for accelerating the convergence of numerical sequences. Each method is based on a systematic approach to reduce or eliminate the leading-order error terms in a convergent or slowly convergent sequence, thereby producing more accurate approximations to the sequence's limit with fewer iterations.
\\
Richardson extrapolation leverages approximations computed with varying discretization steps to cancel dominant error terms, making it particularly effective in numerical differentiation, integration, and the solution of differential equations. Aitken's $\Delta^2$ process provides a simple and computationally efficient approach for accelerating linearly convergent sequences, while the Shanks transformation generalizes this idea by using determinants to improve convergence for more complex series. Wynn's $\epsilon$-algorithm further extends these concepts by providing a recursive scheme capable of handling sequences with oscillatory or nonlinear convergence behavior.
Scalar extrapolation methods play a crucial role in numerical analysis by improving the efficiency and accuracy of iterative computations. 
Vector sequence transformations, including polynomial extrapolation methods such as MPE, RRE, and MMPE, as well as the vector and block $\epsilon$-algorithms, provide powerful tools for accelerating the convergence of sequences arising in iterative numerical methods. These techniques are particularly valuable when dealing with vector sequences generated by the solution of linear or nonlinear systems, where direct methods may be computationally expensive or infeasible. For these class  of vector sequences, Anderson's   method is also a suitable numerical method. We notice also that vector extrapolation methods were generalized to matrix  and tensor cases. The study and implementation of vector sequence transformation techniques are essential components of modern numerical analysis. They provide a rigorous framework for improving convergence, optimizing computational resources, and achieving higher accuracy in a wide range of scientific and engineering applications.

%==========================
\end{document}